\documentclass[11pt]{article}

\usepackage[a4paper,margin=0.72in]{geometry}
\usepackage[T1]{fontenc}
\usepackage[utf8]{inputenc}
\usepackage{amsmath,amssymb,amsthm,bm}
\usepackage{newtxtext,newtxmath}
\usepackage{graphicx}
\usepackage{booktabs}
\usepackage{multirow}
\usepackage{algorithm}
\usepackage{algpseudocode}

\usepackage{xcolor}
\definecolor{linkblue}{RGB}{0,78,155}
\usepackage{hyperref}
\hypersetup{
    colorlinks=true,
    allcolors=linkblue
}
\usepackage{tikz}
\usepackage{pgfplots}
\pgfplotsset{compat=1.18}
\usepackage{caption}
\usepackage{subcaption}
\usetikzlibrary{positioning}
\usepackage{natbib}
\bibpunct[, ]{(}{)}{,}{a}{}{,}

\def\newblock{\ }

\newcommand{\R}{\mathbb{R}}
\providecommand{\argmin}{}
\renewcommand{\argmin}{\operatornamewithlimits{argmin}}
\newcommand{\vmu}{\boldsymbol\mu}
\newcommand{\vone}{\mathbf1}
\newcommand{\rank}{\mathrm{rank}}
\newcommand{\Halmos}{\hfill$\square$}

\newcommand{\RUNAUTHOR}[1]{}
\newcommand{\RUNTITLE}[1]{}
\newcommand{\TITLE}[1]{}
\newcommand{\ARTICLEAUTHORS}[1]{}
\newcommand{\AUTHOR}[1]{}
\newcommand{\AFF}[1]{}
\newcommand{\EMAIL}[1]{}
\newcommand{\ABSTRACT}[1]{}
\newcommand{\FUNDING}[1]{}
\newcommand{\KEYWORDS}[1]{}

\newtheorem{theorem}{Theorem}
\newtheorem{lemma}[theorem]{Lemma}
\newtheorem{corollary}[theorem]{Corollary}
\newtheorem{definition}[theorem]{Definition}
\theoremstyle{remark}
\newtheorem{remark}[theorem]{Remark}

\newcommand{\ArxivTitle}{%
Scalable Mean--Variance Portfolio Optimization\\
via Subspace Embeddings and GPU-Friendly\\
Nesterov-Accelerated Projected Gradient}
\newcommand{\ArxivAuthors}{%
Yi-Shuai Niu\\
Beijing Institute of Mathematical Sciences and Applications (BIMSA)\\
\texttt{niuyishuai@bimsa.cn}\\[0.9em]
Yajuan Wang\\
Beijing Institute of Mathematical Sciences and Applications (BIMSA)\\
\texttt{wangyajuan@bimsa.cn}}
\newcommand{\ArxivDate}{March 29, 2026}

\makeatletter
\renewcommand{\maketitle}{%
  \begin{center}
    {\LARGE\bfseries \ArxivTitle\par}
    \vspace{1em}
    {\normalsize \begin{tabular}{c}\ArxivAuthors\end{tabular}\par}
    \vspace{0.75em}
    {\small \ArxivDate\par}
  \end{center}
  \vspace{0.75em}
  \begin{abstract}
  We develop a sketch-based factor reduction and a Nesterov-accelerated projected
  gradient algorithm (NPGA) with GPU acceleration, yielding a doubly accelerated
  solver for large-scale constrained mean--variance portfolio optimization.
  Starting from the sample covariance factor $L$, the method combines randomized
  subspace embedding, spectral truncation, and ridge stabilization to construct
  an effective factor $L_{\mathrm{eff}}$. It then solves the resulting constrained
  problem with a structured projection computed by scalar dual search and
  GPU-friendly matrix--vector kernels, yielding one computational pipeline for
  the baseline, sketched, and Sketch--Truncate--Ridge (STR)-regularized models.
  We also establish approximation, conditioning, and stability guarantees for the
  sketching and STR models, including explicit $O(\varepsilon)$ bounds for the
  covariance approximation, the optimal value error, and the solution perturbation
  under $(\varepsilon,\delta)$-subspace embeddings. Experiments on synthetic and
  real equity-return data show that the method preserves objective accuracy while
  reducing runtime substantially. On a 5440-asset real-data benchmark with 48374
  training periods, NPGA-GPU solves the unreduced full model in $2.80$ seconds
  versus $64.84$ seconds for Gurobi, while the optimized compressed GPU variants
  remain in the low-single-digit-second regime. These results show that the full
  dense model is already practical on modern GPUs and that, after compression,
  the remaining bottleneck is projection rather than matrix--vector multiplication.
  \end{abstract}
  \noindent\textbf{Keywords.}
  Large-scale mean--variance portfolio optimization, subspace embedding,
  randomized sketching, accelerated projected gradient, GPU computing
  \medskip

  \noindent\textbf{MSC2020 subject classifications.}
  91G10, 90C20, 90C25, 90C06, 65K05, 65F55, 65Y05
  \par\bigskip
}
\makeatother

\newenvironment{APPENDICES}{%
  \appendix
  \setcounter{section}{0}
}{}

\begin{document}

\maketitle

\section{Introduction}
\label{sec:intro}

The mean--variance (MV) portfolio optimization model, first proposed in~\cite{Markowitz1952}, remains a cornerstone of quantitative
investment and risk management.
By minimizing portfolio variance under expected-return and budget constraints,
the MV framework offers a principled balance between risk and reward.
Despite its conceptual simplicity, large-scale implementations in modern markets present both
\emph{statistical} and \emph{computational} challenges.
In today's data-rich environments, the number of investable assets $n$
can reach thousands, while the historical or high-frequency sample length $T$
may span tens or hundreds of thousands, making classical QP methods
memory-intensive and time-consuming.

\paragraph{Data structure and scaling regimes.}
Let $R\in\mathbb{R}^{n\times T}$ be the matrix of asset returns and
$\vmu := \tfrac{1}{T}R\vone_T$ the sample mean vector, where
$\vone_T$ is a $T$-dimensional vector of ones.
Define the centered return (covariance factor) matrix
\[
L := \frac{R-\vmu\vone_T^\top}{\sqrt{T-1}}\in\mathbb{R}^{n\times T},
\qquad
\Sigma = L L^\top.
\]
This formulation highlights the dependence of $\Sigma$ on both
the asset dimension $n$ and the sample length $T$.
When $T$ is comparable to or smaller than $n$, $\Sigma$ is typically ill-conditioned:
small eigenvalues amplify estimation noise and produce unstable weights.
When $T\gg n$, the estimate is statistically more reliable, but the covariance matrix becomes dense and expensive to store and manipulate.
When both $n$ and $T$ are large, these two difficulties appear at the same time.
This is the \emph{statistical--computational double bottleneck} that motivates the present work.

\paragraph{Mathematical formulation and challenges.}
A standard constrained mean--variance portfolio optimization problem is
\begin{equation}
\label{eq:mv}
    \min_{x\in F}\; f(x) := x^\top \Sigma x,
\end{equation}
where
\[
F := \{\,x\in\mathbb{R}^n:\ \vmu^\top x \ge R_{\mathrm{target}},\
       \vone^\top x = 1,\ x \ge 0\,\}
\]
imposes the expected-return, budget, and nonnegativity constraints.
Here, $\vmu$ is the expected return vector, $R_{\mathrm{target}}$ the target
portfolio return, and $x$ the portfolio weights.

From a computational viewpoint, forming $\Sigma=L L^\top$ requires $O(n^2T)$
operations and $O(n^2)$ storage, while each iteration of a QP solver
(e.g., \textsc{Gurobi}, \textsc{MOSEK}, \textsc{CPLEX})
involves dense matrix operations of size $n\times n$,
leading to overall complexity around $O(n^{3.5})$.
Even with modern hardware, these costs are prohibitive for
real-time rebalancing or large-universe simulations.

In short, large-scale MV optimization is difficult for two reasons.
The first is statistical ill-conditioning, caused by limited samples or strong cross-sectional dependence.
The second is computational cost, caused by dense quadratic structure and repeated factorizations in high dimensions.
An effective solver must therefore be both
\emph{statistically robust}, in the sense of being insensitive to sampling noise, and
\emph{computationally scalable}, in the sense of making effective use of modern parallel hardware.

\paragraph{Existing approaches.}
Prior research on scaling mean--variance optimization has evolved along two complementary directions:
(i) \emph{statistically robust covariance estimation}, and
(ii) \emph{computationally scalable optimization algorithms}.
The latter encompasses both classical first-order and distributed methods as well as
modern randomized linear-algebra approaches that compress data before optimization.

\emph{Statistical regularization.}
Since Markowitz~\cite{Markowitz1952}, numerous works have sought to improve
covariance reliability in high dimensions.
Shrinkage estimators toward structured targets provide
well-conditioned estimates~\cite{LedoitWolf2004},
and nonlinear spectral shrinkage further refines eigenvalue dispersion~\cite{LedoitWolf2012}.
Factor and sparsity-based models exploit low-rank market structures~\cite{FanLiaoMincheva2013},
while portfolio constraints can act as implicit regularizers~\cite{JagannathanMa2003}.
Robust formulations explicitly model uncertainty in $\vmu$ and $\Sigma$
to yield stable portfolios under worst-case scenarios~\cite{BenTalNemirovski2000,GoldfarbIyengar2003}.
Empirical studies also find that naive diversification may rival optimized portfolios
under strong estimation noise~\cite{DeMiguelGarlappiUppal2009}.

\emph{First-order and distributed optimization.}
To mitigate computational cost, first-order and splitting methods
(e.g., accelerated gradient, proximal, and ADMM~\cite{Nesterov2004,Bertsekas1999,BeckTeboulle2009,BoydADMM2011})
offer cheap, parallelizable iterations.
Projection routines enforcing simplex or affine constraints are central submodules~\cite{Duchi2008,Condat2016,BauschkeBorwein1996}.
Nevertheless, commercial interior-point solvers such as
Gurobi~\cite{Gurobi}, MOSEK~\cite{MOSEK}, and CPLEX~\cite{CPLEX}
remain memory- and factorization-bound for very large $n$.

\emph{Randomized numerical linear algebra.}
A complementary direction compresses the data before optimization.
Johnson--Lindenstrauss projections~\cite{JohnsonLindenstrauss1984}
and subspace embeddings for regression~\cite{Sarlos2006,ClarksonWoodruff2017}
preserve geometry while drastically reducing dimension.
Surveys in \cite{Woodruff2014,DrineasMahoney2016}
provide guarantees for low-rank approximation and covariance estimation.
These insights motivate constructing a sketched factor
$\tilde L=L\Phi$ such that $\tilde\Sigma=\tilde L\tilde L^\top$
spectrally approximates $\Sigma$, lowering memory and iteration costs
while improving conditioning.

\paragraph{Our approach and contributions.}
This paper develops a computational and theoretical framework for large-scale MV optimization in which approximation, conditioning, and
first-order optimization are designed together rather than treated as
separate steps. The central object is a structured effective factor
$L_{\mathrm{eff}}$ obtained from the sample covariance factor $L$ through
randomized sketching, spectral truncation, and ridge stabilization.
Given $L_{\mathrm{eff}}$, we design a GPU-friendly Nesterov-accelerated
projected gradient method that solves the MV
 problem without
forming the dense covariance matrix explicitly.

The paper makes four main contributions.
First, we introduce a Sketch--Truncate--Ridge (STR) pipeline that compresses the temporal dimension while
explicitly controlling conditioning, which addresses the practical failure of
sketching alone on real financial spectra.
Second, we develop a GPU-friendly Nesterov-accelerated projected gradient algorithm (NPGA) for the baseline, sketched, and STR models, together with factorized
gradient formulas, a convergence analysis for the exact NPGA framework and a discussion of the projector-accuracy issues arising in implementation, and a GPU
implementation path tailored to large dense portfolio instances.
Third, we establish approximation, conditioning, and stability guarantees for the sketching, STR, and factorized models; in particular, under $(\varepsilon,\delta)$-subspace embeddings we derive explicit $O(\varepsilon)$ bounds for the covariance approximation, the optimal value error, and the solution perturbation.
Fourth, we provide empirical evidence on synthetic and real market datasets showing that
the proposed pipeline preserves optimization accuracy, improves numerical
conditioning, and reduces runtime relative to standard dense
formulations.

\paragraph{Paper organization.}
Section~\ref{sec:pipeline} introduces the factor representation,
real-spectrum motivation, STR pipeline, and factorized effective models that
define the computational workflow. Section~\ref{sec:npga} presents the NPGA
solver, its projection routine, convergence statement, and GPU implementation
considerations. Section~\ref{sec:guarantees} then develops the theoretical
analysis for the sketching and STR models, including approximation,
conditioning, and stability guarantees. Section~\ref{sec:experiments} reports numerical results on
synthetic and real datasets, and Section~\ref{sec:conclusion} concludes with
remarks on future directions.

\section{Problem Setup and Computational Pipeline}
\label{sec:pipeline}

This section presents the computational workflow used throughout the paper.
We start from the factor representation $\Sigma = LL^\top$, explain why sketching alone is insufficient on real financial spectra, define the STR approximation that produces a structured effective factor $L_{\mathrm{eff}}$, and show how that factorization leads directly to a scalable first-order solver.
The supporting theoretical analysis is deferred to Section~\ref{sec:guarantees}.

The workflow is simple.
We represent the sample covariance through the tall factor $L$, compress the temporal dimension by randomized embedding to obtain $\tilde L = L\Phi$, truncate noisy singular modes and add ridge regularization to obtain a numerically stable factor $L_{\mathrm{eff}}$, and finally solve the constrained MV problem by NPGA using only matrix--vector products with $L_{\mathrm{eff}}$ and $L_{\mathrm{eff}}^\top$.

\subsection{Sample Covariance and Factorization}
Let $R\in\mathbb{R}^{n\times T}$ be the panel of historical (demeaned later) returns, and let
\[
\bar r \;=\;\frac{1}{T}\sum_{t=1}^T R(:,t)\;\in\;\mathbb{R}^n,\qquad
R_c \;=\; R - \bar r\,\vone_T^\top \;\in\;\mathbb{R}^{n\times T}.
\]
Define
\[
L \;=\; \frac{1}{\sqrt{T-1}}\,R_c \;\in\;\mathbb{R}^{n\times T},
\]
so that the (unbiased) sample covariance is represented exactly as
\[
\Sigma \;=\; \frac{1}{T-1}\,R_cR_c^\top \;=\; L\,L^\top.
\]
Since centering removes at most one degree of freedom, we have
\[
\rank(\Sigma)=\rank(LL^\top)\le \rank(L)\le \min\{n,T-1\}.
\]
When $T$ is large, explicitly storing $L$ or forming $\Sigma$ is expensive:
$O(nT)$ memory for $L$, or $O(n^2)$ for $\Sigma$, and $O(n^2T)$ time to form
$LL^\top$.

These observations motivate the need for structured approximations of the
sample covariance that preserve the essential spectral information while
reducing storage and computational cost. In the next subsection, we examine
the spectrum of a real-world covariance matrix constructed from large-scale
market data, which will guide our subsequent algorithmic design.

\subsection{Spectral Characteristics of Real-World Covariance Matrices}
\label{sec:real-spectra}

We now examine the spectral
properties of a large real-world covariance matrix of the form
$\Sigma = LL^\top$. In our empirical study, the raw intraday returns come from
the RESSET database and cover \emph{5-minute} returns for A--share stocks
listed on the Shanghai and Shenzhen exchanges from
\emph{January~2,~2020 to June~30,~2025}. The vendor files are aligned to a
common 5-minute trading grid, and missing return entries are filled with zero
when the balanced panel is assembled. After removing assets with insufficient
trading histories, the resulting cleaned matrix contains $(n,T)=(5440,66412)$
observations, where $T$ is the total number of 5-minute intervals over the
sample horizon.

\begin{figure}[t]
  \centering
  \includegraphics[width=0.52\linewidth]{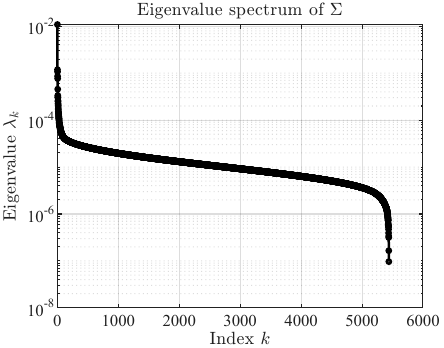}
  \caption{Eigenvalue spectrum of the real-world covariance matrix.
           The spectrum decays rapidly for the first few dozen components
           before entering a long, slowly decaying bulk spanning several
           thousand dimensions.}
  \label{fig:spectrum-real}
\end{figure}

Let $\Sigma = U \Lambda U^\top$ be the eigen-decomposition, where
$\Lambda = \mathrm{diag}(\lambda_1,\dots,\lambda_n)$ with
$\lambda_1\ge\cdots\ge\lambda_n\ge0$.
Figure~\ref{fig:spectrum-real} displays the eigenvalue spectrum of ~$\Sigma$ on a logarithmic scale.
The leading eigenvalues exhibit a sharp drop, corresponding to
a small number of dominant market and industry factors.
However, after the first few dozen components, the spectrum enters
a \emph{very long and nearly flat bulk region} extending over several
thousand eigenvalues.
A geometric fit over this region yields an average decay factor of $0.9995$,
indicating almost linearly in the log domain. Near the tail end of the spectrum, the eigenvalues collapse rapidly toward 
the numerical noise floor.  
This structure is characteristic of high-dimensional financial datasets:
\emph{a handful of strong factors followed by a massive high-dimensional
noise component.}

\begin{figure}[t]
  \centering
  \includegraphics[width=0.52\linewidth]{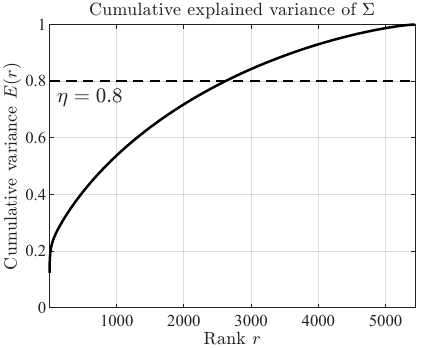}
  \caption{Cumulative explained variance $E(r)$.
           Due to the extremely long noisy bulk, $80\%$ of the total
           variance is reached only after more than $2500$ eigenvalues,
           demonstrating the inadequacy of energy-based rank selection
           for real-world covariance matrices.}
  \label{fig:energy-real}
\end{figure}

\paragraph{Failure of energy-based rank truncation.}
A common approach for defining an effective rank is the cumulative
explained variance
\[
  E(r)
  =
  \frac{\sum_{k=1}^r \lambda_k}{\sum_{k=1}^n \lambda_k}.
\]
However, the bulk region contributes substantial total energy due to
its large dimensionality.
For the present dataset, even a moderate threshold $\eta=0.8$
requires retaining
\[
  r_\eta = \min\{\, r : E(r)\ge\eta \,\}
  =
  2616,
\]
which accounts for nearly half the ambient dimension (see Figure~\ref{fig:energy-real}).
Such ranks are far too large to yield computational savings and
do not meaningfully separate structural factors from high-dimensional
noise.

\paragraph{Implications for numerical conditioning.}
The slow spectral decay of $\Sigma$ results in severe ill-conditioning.
Even when restricted to the first $r_\eta=2616$ eigen-directions,
the effective condition number remains
\[
    \kappa_{\mathrm{eff}}
    =
    \frac{\lambda_1}{\lambda_{r_\eta}}
    \approx 10^3,
\]
while the full covariance matrix of order $10^5$ is substantially more ill-conditioned.
Consequently, gradient-based optimization algorithms that rely on
$\Sigma x$ suffer from unstable search directions and slow convergence
when applied directly to~$\Sigma$ or to high-rank approximations
obtained via energy truncation.

\paragraph{Motivation for structured approximation.}
The empirical spectrum reveals three distinct regions: a small number of significant factors at the head, an extremely long and slowly decaying noisy bulk, and a numerically unstable tail of near-zero eigenvalues.
This structure cannot be captured by naive low-rank truncation
or explained-variance criteria.
Instead, an effective approximation must
preserve the structural factor subspace,
compress the bulk into a stable isotropic component,
and regularize the spectral tail.
These observations motivate the Sketch--Truncate--Ridge (STR) procedure
introduced in the subsequent sections.

\subsection{Randomized Subspace Embedding}
\label{sec:embedding}

As shown in subsection~\ref{sec:real-spectra}, the spectrum of 
$\Sigma = LL^\top$ contains only a few informative directions, while the 
remaining thousands fall into a slowly decaying noisy bulk.  
To extract the factor subspace efficiently without explicitly forming 
$L$ or $\Sigma$, we employ \emph{randomized subspace embedding}, which 
reduces the temporal dimension from $T$ to $s\ll T$ while preserving all 
quadratic forms $x^\top L L^\top x$. This reduction serves as the 
computational backbone for all subsequent approximation stages.

To obtain such a compressed representation, we post-multiply 
$L\in\mathbb{R}^{n\times T}$ by a random sketching matrix 
$\Phi\in\mathbb{R}^{T\times s}$, yielding
\[
\tilde L = L\Phi \in \mathbb{R}^{n\times s},
\qquad
\tilde\Sigma = \tilde L \tilde L^\top \approx \Sigma.
\]
The sketching matrix $\Phi$ is chosen so that, with probability at least 
$1-\delta$, it acts as an $(\varepsilon,\delta)$-subspace embedding for 
$\mathrm{Im}(L^\top)$:
\[
(1-\varepsilon)\,\|L^\top x\|_2^2
\;\le\;
\|\tilde L^\top x\|_2^2
\;\le\;
(1+\varepsilon)\,\|L^\top x\|_2^2,
\qquad \forall\,x\in\mathbb{R}^n.
\]
Thus, the Euclidean geometry of $L^\top$ and hence all quadratic forms 
$x^\top\Sigma x$ is preserved up to a multiplicative distortion 
of $(1\pm\varepsilon)$.

We employ two canonical constructions of $\Phi$ from
randomized numerical linear algebra:

\paragraph{(1) Gaussian Johnson--Lindenstrauss (JL) projection \cite{JohnsonLindenstrauss1984}.}
Each entry of $\Phi$ is drawn i.i.d.\ as
$\Phi_{ij}\sim\mathcal{N}(0,1/s)$.
This dense random projection requires $O(Ts)$ time and $O(Ts)$ storage,
and satisfies the Johnson--Lindenstrauss lemma:
for any fixed finite set of vectors,
Euclidean distances are preserved within a factor of $(1\pm\varepsilon)$
with probability at least $1-\delta$
when
\[
s = O\!\bigl(\varepsilon^{-2}\log(1/\delta)\bigr).
\]
In practice, Gaussian JL sketches are simple to generate and yield
tight theoretical error bounds, but may be computationally expensive
for extremely large $T$ due to their dense structure.

\paragraph{(2) Count--Sketch mapping \cite{Sarlos2006,Woodruff2014,ClarksonWoodruff2017}.}
Count--Sketch offers a sparse and highly scalable alternative.
Each row of $\Phi$ contains exactly one nonzero entry $\pm1$
at a uniformly hashed column position:
\[
\Phi(i,\,h(i)) = s_i,\quad i=1,\ldots,T,
\]
where $h:\{1,\ldots,T\}\to\{1,\ldots,s\}$ is a random hash function
and $s_i\in\{\pm1\}$ are independent Rademacher signs.
This structure yields only $O(T)$ time and $O(T)$ storage
to generate and apply $\Phi$, making it well suited for large-scale
or streaming data.
Count--Sketch achieves an $(\varepsilon,\delta)$-subspace embedding
for any $r$-dimensional subspace with
\[
s = O\!\bigl((r+\log(1/\delta))/\varepsilon^2\bigr),
\]
up to modest constant factors compared with Gaussian JL.
While the asymptotic dependence is the same, the sparse structure
dramatically reduces memory footprint and multiplication cost.

\paragraph{Discussion.}
Both constructions guarantee that the resulting sketched matrix
$\tilde L=L\Phi$ preserves the spectral geometry of $L$ with
probability at least $1-\delta$.
Gaussian JL offers optimal theoretical distortion,
whereas Count--Sketch trades a slightly looser constant for massive gains
in efficiency and scalability.
Throughout this paper, either choice of $\Phi$ is admissible,
and all subsequent results on spectral approximation, robustness,
and conditioning apply uniformly to both types of embeddings.
The next subsection formalizes these properties through the
Johnson--Lindenstrauss lemma and the
$(\varepsilon,\delta)$-subspace embedding theorem. These guarantees are essential because the STR procedure that follows
relies on $\tilde L$ as a surrogate for $L$ in both spectral truncation
and regularization.  
We therefore summarize the foundational results below.

The empirical spectral evidence above motivates a conditioning mechanism beyond pure sketching. We next introduce the STR pipeline, which combines randomized embedding, spectral truncation, and ridge stabilization to produce a numerically stable effective factor for optimization.

\subsection{STR: Overview and Motivation}
The empirical study in Section~\ref{sec:real-spectra} shows that real-world
covariance matrices exhibit a highly unbalanced spectrum:
a handful of dominant factors, followed by an extremely long,
slowly decaying bulk and a numerically unstable tail.
Such spectral geometry leads to severe ill-conditioning and renders
gradient-based solvers inefficient or unstable.

Randomized subspace embedding (Section~\ref{sec:embedding}) provides an
efficient way to compress the temporal dimension while preserving all
quadratic forms $x^\top\Sigma x$ up to a factor of $(1\pm\varepsilon)$.
However, sketching alone does \emph{not} improve conditioning.
Indeed, from the positive semidefinite (PSD) ordering
\[
(1-\varepsilon)\Sigma \preceq \tilde\Sigma \preceq (1+\varepsilon)\Sigma,
\]
one obtains
\[
(1-\varepsilon)\,\lambda_{\min}(\Sigma)\ \le\ \lambda_{\min}(\tilde\Sigma)\ \le\ \lambda_{\max}(\tilde\Sigma)\ \le\ (1+\varepsilon)\,\lambda_{\max}(\Sigma).
\]
Consequently, the condition number satisfies
\[
\kappa(\tilde\Sigma)
\;=\;\frac{\lambda_{\max}(\tilde\Sigma)}{\lambda_{\min}(\tilde\Sigma)}
\;\le\;
\frac{1+\varepsilon}{1-\varepsilon}\,\kappa(\Sigma),
\]
which bounds the possible distortion but does not guarantee any reduction in
$\kappa(\Sigma)$. Thus the small singular modes of $L$, which generate the long noisy bulk
of $\Sigma$, continue to dominate numerical performance even after sketching.

These observations motivate a structured procedure that not only preserves
the geometry of $LL^\top$ but also reshapes its spectrum into a numerically
stable form. Since the eigenvalues of $\Sigma$ are given by the squared
singular values of $L$,
\[
\lambda_i(\Sigma)=\sigma_i(L)^2,\qquad i=1,\dots,\rank(L),
\]
the ill-conditioning of $\Sigma$ originates directly from the singular-value
spread of $L$.  Thus, improving the conditioning of $\Sigma$ requires
manipulating the singular spectrum of $L$ itself.

\subsection{STR: Description of the Procedure}
To this end, we introduce a three-stage conditioning framework,
\emph{Sketch--Truncate--Ridge (STR)}, which produces a well-conditioned and
structurally faithful approximation of $\Sigma$.

\paragraph{(1) Sketch stage.}
Apply a Johnson--Lindenstrauss transform or a CountSketch matrix
$\Phi\in\mathbb{R}^{T\times s}$ ($s\ll T$) and form
\[
\tilde L = L\Phi,\qquad \tilde\Sigma=\tilde L\tilde L^\top.
\]
This step compresses the temporal dimension and provides a
$(1\pm\varepsilon)$ spectral approximation to $\Sigma$, but does not by itself
improve conditioning.

\paragraph{(2) Truncation stage.}
Compute the thin singular value decomposition (SVD)
\[
\tilde L=USV^\top,\qquad
S=\mathrm{diag}(\sigma_1,\dots,\sigma_r),\quad
\sigma_1\ge\cdots\ge\sigma_r>0.
\]
Let $U_\ell,S_\ell,V_\ell$ denote the first $\ell$ singular components and
construct the rank-$\ell$ approximation
\[
\tilde L_\ell := U_\ell S_\ell V_\ell^\top .
\]
The truncation step removes small singular modes of $\tilde L$ and raises
the smallest retained singular value from $\sigma_r$ to $\sigma_\ell$,
thereby eliminating directions that would otherwise produce near-zero
eigenvalues in $\Sigma$. In practice, $\ell$ is chosen at the end of the
dominant spectral core of $\tilde L$ rather than by a cumulative-variance
rule; a detailed selection rule and a computable truncation bound are given
in Appendix~\ref{app:str-discussion}.

\paragraph{(3) Ridge stage.}
To restore full rank and ensure numerical stability, we apply a uniform
spectral lift:
\[
\widehat\Sigma := \tilde L_\ell \tilde L_\ell^\top + \gamma I_n,\qquad \gamma>0.
\]
The eigenvalues of $\widehat\Sigma$ are
\[
\{\sigma_1^2+\gamma,\dots,\sigma_\ell^2+\gamma,
\underbrace{\gamma,\dots,\gamma}_{n-\ell}\},
\]
so that
\[
\lambda_{\min}(\widehat\Sigma)=\gamma,
\qquad
\kappa(\widehat\Sigma)=\frac{\sigma_1^2+\gamma}{\gamma}.
\]
The ridge term fills the zero spectrum created by truncation and guarantees
uniform conditioning.

Finally, note that $\sigma_1^2(\tilde L)$ is controlled by the embedding
property:
\[
\sigma_1^2(\tilde L)
=\lambda_{\max}(\tilde\Sigma)
\le (1+\varepsilon)\lambda_{\max}(\Sigma).
\]
Thus, ensuring $\kappa(\widehat\Sigma)<\kappa(\Sigma)$ requires choosing
\[\boxed{
\gamma \;>\;
\frac{(1+\varepsilon)\lambda_{\min}(\Sigma)\,\lambda_{\max}(\Sigma)}
{\lambda_{\max}(\Sigma)-(1+\varepsilon)\lambda_{\min}(\Sigma)},}
\]
which guarantees that STR strictly improves the conditioning of $\Sigma$.

The STR procedure transforms the raw return matrix $R$ into a stable and
well-conditioned covariance approximation through sketching, truncation, and
ridge stabilization. A schematic summary of this pipeline is provided in
Figure~\ref{fig:app-str-pipeline} in the appendix.
The resulting approximation enters the optimization model only through factorized matrix--vector products. This viewpoint leads to reduced effective models and cheap gradients, which will be exploited directly by the solver in Section~\ref{sec:npga}.

The quadratic objective $f(x)=x^\top\Sigma x$ in the MV problem~\eqref{eq:mv} involves
the dense covariance matrix $\Sigma\in\mathbb{R}^{n\times n}$.
Directly manipulating $\Sigma$ incurs $O(n^2)$ storage and $O(n^3)$
factorization cost. To expose the low-rank structure of $L$, we introduce
an auxiliary variable $z=L^\top x\in\mathbb{R}^{T}$ (or $\tilde z=\tilde L^\top x \in \R^s$
after sketching, $\hat z=\tilde L_\ell^\top x \in \R^\ell$
after STR). The MV problem can then be
reformulated as
\begin{equation}
\label{eq:Lz-problem}
\min
\;\; \|z\|_2^2
\quad
\text{s.t.}\quad
L^\top x=z, z\in\mathbb{R}^{T}, x\in F,
\end{equation}
where $F=\{x\in\mathbb{R}^n:\vmu^\top x\ge R_{\mathrm{target}},
\vone^\top x=1,\,x\ge0\}$ is the feasible set. This formulation replaces the dense quadratic term by a least-norm objective
and an explicit linear coupling constraint. The problem remains convex
because both $\|z\|_2^2$ and the affine constraint are convex, and the feasible
region of $(x,z)$ is polyhedral.

The $L$--$z$ reformulation offers two main advantages:
it enables dimensional separation between the asset and temporal spaces,
and it preserves a factorized structure that allows
efficient gradient computation under sketching and STR.

\subsection{Dimensional separation and sketched models}
  The coupling $L^\top x = z$ explicitly separates the
  \emph{asset space} ($x\in\R^n$) and the
  \emph{temporal space} ($z\in\R^T$),
  allowing subspace embedding to be applied \emph{only} on the temporal
  dimension.

  \smallskip
  \noindent\emph{(a) Sketched model.}
  After applying a random projection $\Phi\in\R^{T\times s}$ with $s\ll T$,
  we obtain the reduced formulation
  \begin{equation}
  \label{eq:sketched-model}
  \min_{x\in F,\;\tilde z\in\R^{s}}
  \|\tilde z\|_2^2
  \quad\text{s.t.}\quad
  \tilde L^\top x = \tilde z,
  \qquad \tilde L := L\Phi.
  \end{equation}
  Eliminating $\tilde z$ gives the sketched quadratic form
  \[
  \min_{x\in F}\; \tilde f(x):=x^\top\tilde\Sigma x = \|\tilde L^\top x\|_2^2,
  \qquad
  \tilde\Sigma := \tilde L\tilde L^\top.
  \]

  \smallskip
  \noindent\emph{(b) STR model.}
  If the STR procedure is applied, then
  \begin{equation}
  \label{eq:STR-model}
  \min_{x\in F,\;\hat z\in\R^{k}}
  \|\hat z\|_2^2 + \gamma\|x\|_2^2
  \quad\text{s.t.}\quad
  \tilde L_\ell^\top x = \hat z,
  \end{equation}
  where $\tilde L_\ell\in\R^{n\times \ell}$ is the truncated and regularized
  factor obtained from STR.
  Eliminating $\hat z$ yields the regularized quadratic form
  \[
  \min_{x\in F}\; \widehat f(x):=x^\top\widehat\Sigma x = \|\tilde L_\ell^\top x\|_2^2 + \gamma \|x\|_2^2,
  \qquad
  \widehat\Sigma := \tilde L_\ell\tilde L_\ell^\top + \gamma I_n.
  \]
  Hence STR reduces the temporal dimension from $T$ to $\ell$
  and stabilizes the covariance spectrum through ridge lifting.

\subsection{Structure-exploiting gradient}
  The factorized representation enables fast gradient computation
  for all three cases.

  \smallskip
  \noindent\emph{(a) Baseline (no sketch).}
  For $f(x)=x^\top\Sigma x$ with $\Sigma=LL^\top$,
  \[
  \nabla f(x) = 2LL^\top x = 2Lz,
  \qquad z = L^\top x,
  \]
  requiring two matrix--vector products with $L$ and $L^\top$,
  each of cost $O(nT)$.

  \smallskip
  \noindent\emph{(b) Sketched gradient.}
  For the sketched objective $\tilde f(x)=x^\top\tilde\Sigma x$ with
  $\tilde\Sigma=\tilde L\tilde L^\top$,
  \[
  \nabla \tilde f(x)
  = 2\tilde L\tilde L^\top x
  = 2\tilde L\tilde z,
  \qquad \tilde z=\tilde L^\top x.
  \]
  This reduces the gradient cost from $O(nT)$ to $O(ns)$
  since $s\ll T$.

  \smallskip
  \noindent\emph{(c) STR gradient.}
  For the STR-regularized objective
  $\widehat f(x)=x^\top\widehat\Sigma x$
  with $\widehat\Sigma=\tilde L_\ell\tilde L_\ell^\top+\gamma I_n$,
  \[
  \nabla \widehat f(x)
  = 2\tilde L_\ell(\tilde L_\ell^\top x) + 2\gamma x
  = 2\tilde L_\ell\hat z + 2\gamma x,
  \qquad \hat z = \tilde L_\ell^\top x.
  \]
  The cost per iteration is $O(n\ell)$, dominated by two
  matrix--vector multiplications with $\tilde L_\ell$ and $\tilde L_\ell^\top$.
  Compared with the baseline $O(nT)$ and sketched $O(ns)$,
  STR achieves the lowest complexity and improved conditioning.

\section{GPU-Friendly Nesterov-Accelerated Projected Gradient Algorithm (NPGA)}
\label{sec:npga}
\subsection{NPGA: A Unified Scheme for Baseline, Sketch, and STR}

We parameterize all three models (the $L$--$z$ equivalent formulation~\eqref{eq:Lz-problem}, its sketched-model \eqref{eq:sketched-model} and STR-model \eqref{eq:STR-model}) by an effective factor $L_{\mathrm{eff}}$ and a ridge $\gamma\ge 0$:
\[
L_{\mathrm{eff}}\in\{L,\ \tilde L,\ \tilde L_\ell\},\qquad
\gamma\in\{0,\ 0,\ >0\}.
\]
The unified objective and gradient read
\begin{equation}\label{eq:npga-unified-grad}
f(x)=\|L_{\mathrm{eff}}^\top x\|_2^2+\gamma\|x\|_2^2,
\qquad
\nabla f(x)=2\,L_{\mathrm{eff}}\!\left(L_{\mathrm{eff}}^\top x\right)+2\gamma x.
\end{equation}

\begin{table}[ht!]
\centering
\caption{Per-iteration complexity of NPGA. }
\label{tab:model-mapping}
\begin{tabular}{llcc}
\toprule
\textbf{Model} & \textbf{$L_{\mathrm{eff}}$} & \textbf{$\gamma$} & \textbf{Cost / iter} \\
\midrule
Baseline \eqref{eq:Lz-problem} & $L\in\R^{n\times T}$ & $0$ & $O(nT)$ \\
Sketch \eqref{eq:sketched-model}  & $\tilde L\in\R^{n\times s}$ & $0$ & $O(ns)$ \\
STR \eqref{eq:STR-model}   & $\tilde L_\ell\in\R^{n\times \ell}$ & $>0$ & $O(n\ell)$ \\
\bottomrule
\end{tabular}
\end{table}

We now propose an efficient first-order algorithm, the GPU-Friendly
Nesterov-Accelerated Projected Gradient Algorithm (NPGA), that fully exploits
the problem structure revealed by the $L$--$z$ formulation. The proposed NPGA
reads as in Algorithm~\ref{alg:npga}. Let $\alpha_k$ denote the step size at
iteration $k$ and $y^k$ the extrapolated Nesterov point. The step sizes may be
chosen either by a fixed Lipschitz-based rule or by an adaptive backtracking
line search. In the general convex regime we use the standard FISTA-type
momentum update, whereas in the strongly convex regime the same framework also
admits the classical constant-momentum specialization discussed below. We
terminate the main loop when the projected gradient residual satisfies
\[
\|\Pi_F(x^k-\alpha_k\nabla f(x^k)) - x^k\|_2 \le \epsilon.
\]
The final iterate $x^k$ is then returned as the solution.
By leveraging the factorized gradient representation~\eqref{eq:npga-unified-grad},
each iteration involves only two matrix--vector products with
$L_{\mathrm{eff}}$ and $L_{\mathrm{eff}}^\top$ (plus an $O(n)$ ridge term if present),
operations that are fully parallelizable on modern GPUs.
This structural efficiency enables the method to scale to
large, high-dimensional portfolio problems where the covariance matrix
is ill-conditioned and costly to form or store explicitly.
To further accelerate convergence, we incorporate Nesterov's extrapolation,
which yields the standard $O(1/k^2)$ guarantee for smooth convex objectives and
also admits the usual linear-rate specialization in the strongly convex regime,
while maintaining full GPU parallelization efficiency.

Note that the NPGA framework can also be applied to the classical
mean--variance model~\eqref{eq:mv} with an explicit covariance matrix
$\Sigma$ in place of $L_{\mathrm{eff}}$.
In this case, the gradient is $\nabla f(x)=2\Sigma x$, and each iteration
requires $O(n^2)$ operations.
Hence, NPGA serves as a unified GPU-friendly solver for the full spectrum of
mean--variance formulations.

\begin{algorithm}[t]
\caption{NPGA: Nesterov--Accelerated Projected Gradient}
\label{alg:npga}
\begin{algorithmic}[1]
\Require $L_{\mathrm{eff}}$, ridge $\gamma\ge 0$, initial $x^0\in F$, stepsizes $\{\alpha_k\}$, tolerance $\epsilon>0$.
\Ensure $x^k$.
\State \textbf{Initialize:} $y^0=x^0$, $t_0=1$, $k=0$.
\While{ not converged }
  \State Compute gradient: $\nabla f(y^k)=2\,L_{\mathrm{eff}}(L_{\mathrm{eff}}^\top y^k)+2\gamma y^k$.
  \State Gradient step: $v^k=y^k-\alpha_k\nabla f(y^k)$.
  \State Project onto feasible set: $x^{k+1}=\Pi_F(v^k)$.
  \State Update momentum:
        $t_{k+1}=\tfrac{1+\sqrt{1+4t_k^2}}{2}$,
        $\beta_k=\tfrac{t_k-1}{t_{k+1}}$.
  \State Extrapolation: $y^{k+1}=x^{k+1}+\beta_k(x^{k+1}-x^k)$.
  \State $k\leftarrow k+1$.
\EndWhile
\end{algorithmic}
\end{algorithm}

\paragraph{Stepsize policy.}
Let $L_f$ be the Lipschitz constant of $\nabla f$:
\[
L_f=2\,\|L_{\mathrm{eff}}L_{\mathrm{eff}}^\top\|_2+2\gamma
=2\,\|L_{\mathrm{eff}}\|_2^2+2\gamma.
\]
The NPGA framework allows either a fixed-step policy or an adaptive
backtracking policy, both expressed through the same model-dependent
smoothness constant \(L_f\). Algorithm~\ref{alg:npga} displays the default
\(t_k\)-based FISTA update used in the convex regime and in the reported
benchmarks. When a lower curvature bound \(m_f>0\) is available, the same
projected-gradient framework also admits the classical constant-momentum
specialization discussed below.

\paragraph{Fixed step.}
One may use a constant step size \(\alpha_k\equiv \alpha\) with
\(0<\alpha\le 1/L_f\). In practice, \(\|L_{\mathrm{eff}}\|_2\) is estimated
either by the largest singular value or via a few iterations of the power
method, which is GPU-friendly. When a certified admissible step is needed for
the theorem, one may instead use any provable upper bound on \(L_f\) or invoke
the backtracking rule below.

\paragraph{Adaptive backtracking.}
Alternatively, NPGA may use the standard accelerated proximal-gradient
backtracking rule. Given an initial trial step \(\alpha_k^{(0)}>0\) and a
shrinkage factor \(\rho\in(0,1)\), we decrease \(\alpha_k\) geometrically
until the accepted point \(x^{k+1}=\Pi_F(y^k-\alpha_k\nabla f(y^k))\) satisfies
\[
f(x^{k+1})
\le
f(y^k)
+
\langle \nabla f(y^k),x^{k+1}-y^k\rangle
+
\frac{1}{2\alpha_k}\|x^{k+1}-y^k\|_2^2.
\]
By \(L_f\)-smoothness, this condition is guaranteed whenever \(\alpha_k\le
1/L_f\), so the accepted step sizes remain admissible. This gives a flexible
variable-step implementation without changing the basic NPGA framework.

\paragraph{Strongly convex specialization.}
If a lower curvature bound \(m_f>0\) is also available, then the same NPGA
framework admits the standard constant-momentum choice
\[
\beta_k\equiv
\beta:=
\frac{1-\sqrt{\alpha m_f}}
{1+\sqrt{\alpha m_f}},
\]
associated with any fixed admissible step \(0<\alpha\le 1/L_f\). This is the
version that yields the linear-rate statement recorded in
Theorem~\ref{thm:npga-convergence}(b); Algorithm~\ref{alg:npga} itself shows
the default \(t_k\)-based update used in part~(a).

\begin{remark}[Power method for estimating $\|L_{\mathrm{eff}}\|_2$]
The spectral norm $\|L_{\mathrm{eff}}\|_2$ equals the largest singular value
of $L_{\mathrm{eff}}$, which can be efficiently estimated by the classical
\emph{power method} through the iterative update
\[
u^{t+1} = 
\frac{L_{\mathrm{eff}}\!\left(L_{\mathrm{eff}}^\top u^t\right)}
     {\bigl\|L_{\mathrm{eff}}\!\left(L_{\mathrm{eff}}^\top u^t\right)\bigr\|_2},
\qquad 
u^0 = \frac{r}{\|r\|_2},\quad r\sim\mathcal{N}(0,I_n).
\]
The associated Rayleigh quotient
$\|L_{\mathrm{eff}}^\top u^t\|_2$
monotonically increases and converges to $\|L_{\mathrm{eff}}\|_2$
when $L_{\mathrm{eff}}$ has a nonzero spectral gap.
Each iteration involves two matrix--vector products, whose cost per iteration is summarized in Table~\ref{tab:model-mapping} and is readily parallelizable on GPUs.
Empirically, 5--10 iterations are sufficient to obtain a reliable spectral
proxy. In the fixed-step implementation this yields an inexpensive practical
step-size estimate, while the theorem itself may instead rely on any certified
upper bound or on the backtracking rule.
\end{remark}

\paragraph{Projection onto \texorpdfstring{$F$}{F}.}
The feasible set
\[
F=\{x\in\mathbb{R}^n:\ \vmu^\top x\ge R_{\mathrm{target}},\
\vone^\top x=1,\ x\ge0\}
\]
can be written as the intersection
\[
F=\Delta\cap H,
\qquad
\Delta=\{x\in\mathbb{R}^n:\ x\ge0,\ \vone^\top x=1\},
\qquad
H=\{x\in\mathbb{R}^n:\ \vmu^\top x\ge R_{\mathrm{target}}\}.
\]
We compute $\Pi_F(v)$ by exploiting the fact that $F$ is a simplex intersected
with a single return halfspace. Let
$
x_\Delta=\Pi_\Delta(v).
$
If $\vmu^\top x_\Delta\ge R_{\mathrm{target}}$, then $x_\Delta$ is already
feasible for $F$ and therefore coincides with $\Pi_F(v)$. The nontrivial case
is when $x_\Delta$ violates the return constraint, in which case the halfspace
constraint is active at the projector.

\paragraph{Projection onto the simplex.}
We first project onto
$
\Delta.
$
The simplex projection admits the closed form
\citep{Duchi2008,Condat2016}:
find the threshold $\tau$ satisfying
\begin{equation}
\label{eq:simplex-tau}
\sum_{i=1}^n \max\{v_i-\tau,0\}=1,
\end{equation}
and set
\begin{equation}
\label{eq:simplex-proj}
(\Pi_\Delta (v))_i = \max\{v_i-\tau,0\},\quad i=1,\ldots,n.
\end{equation}
The threshold $\tau$ can be obtained by sorting $v$
in descending order and identifying the largest index
$j$ such that
\[
v_{j} - \frac{1}{j}\Big(\sum_{i=1}^j v_{i} - 1\Big) > 0,
\]
then
\(\displaystyle
\tau = \frac{1}{j}\Big(\sum_{i=1}^j v_{i} - 1\Big).
\)
This step costs $O(n\log n)$ because of sorting.

\paragraph{Projection onto the return halfspace.}
The halfspace projection has the closed form
\begin{equation}
\label{eq:proj-halfspace}
\Pi_H(y)=
\begin{cases}
y, & \text{if }\vmu^\top y\ge R_{\mathrm{target}},\\[3pt]
y+\dfrac{R_{\mathrm{target}}-\vmu^\top y}{\|\vmu\|_2^2}\,\vmu,
& \text{otherwise.}
\end{cases}
\end{equation}
This update costs $O(n)$.

\paragraph{Exact projection onto $F$ via a scalar dual search.}
When $\vmu^\top x_\Delta<R_{\mathrm{target}}$, the KKT conditions imply that the
projector lies on the active face $\vmu^\top x=R_{\mathrm{target}}$. Introducing
a Lagrange multiplier $\nu\ge 0$ for the halfspace constraint gives the
one-parameter family
\begin{equation}
\label{eq:projF-dualpath}
x(\nu)
=
\argmin_{x\in\Delta}
\left\{
\frac{1}{2}\|x-v\|_2^2 + \nu\bigl(R_{\mathrm{target}}-\vmu^\top x\bigr)
\right\}
=
\Pi_\Delta(v+\nu \vmu).
\end{equation}
Define
\[
\phi(\nu):=\vmu^\top x(\nu)=\vmu^\top \Pi_\Delta(v+\nu\vmu),
\qquad \nu\ge 0.
\]
The map $\phi$ is continuous and nondecreasing. Hence, whenever
$R_{\mathrm{target}}\le \max_i \mu_i$, there exists $\nu^\star\ge 0$ such that
\[
\vmu^\top x(\nu^\star)=R_{\mathrm{target}},
\qquad
\Pi_F(v)=x(\nu^\star).
\]
Our implementation first brackets $\nu^\star$ by geometric expansion
$\nu=1,2,4,\dots$ until $\phi(\nu)\ge R_{\mathrm{target}}$, and then applies
bisection on the bracket. The iteration terminates once
\[
\bigl|\phi(\nu)-R_{\mathrm{target}}\bigr|
\le
\tau_{\mathrm{proj}}\max\{1,|R_{\mathrm{target}}|\}.
\]
Each evaluation of $\phi(\nu)$ requires one simplex projection and one dot
product, so the dominant cost remains the $O(n\log n)$ simplex sort. After
bracketing, the total work is $O(n\log n\log(1/\tau_{\mathrm{proj}}))$, while
the projection itself is exact up to the scalar-search tolerance.

As a safeguarded fallback, we retain a Dykstra-style alternating-projection method
\citep{BauschkeBorwein1996} that combines $\Pi_\Delta$ and $\Pi_H$ if the
bracketing phase does not succeed within the prescribed iteration budget. In
our experiments this safeguard is rarely triggered; the default path is the
structured scalar-search projector above.
\paragraph{Stopping condition.} We stop when the projected-gradient residual
\[\|\Pi_F(x^k-\alpha\nabla f(x^k)) - x^k\|_2 \le \epsilon,\] which is
equivalent to requiring that the current iterate \(x^k\) approximately
satisfies the KKT conditions. For any fixed
stepsize \(\alpha>0\), define the \emph{projected gradient mapping}
\(G_\alpha(x)=\tfrac{1}{\alpha}\bigl(x-\Pi_F(x-\alpha\nabla f(x))\bigr)\).
In convex optimization, $x^\star$ is optimal if and only if
$0\in\nabla f(x^\star)+N_F(x^\star)$, where $N_F(x^\star)$ is the normal cone to $F$ at $x^\star$.
This condition can be rewritten as $\Pi_F(x^\star-\alpha\nabla f(x^\star))=x^\star$
for any $\alpha>0$, i.e., $G_\alpha(x^\star)=0$.
Hence $\|G_\alpha(x)\|_2$ measures the violation of the KKT stationarity condition,
and vanishes exactly at an optimal point.
The practical criterion
$\|\Pi_F(x^k-\alpha\nabla f(x^k)) - x^k\|_2 \le \epsilon$
thus ensures that both the gradient balance and the projected feasibility residuals are
below the prescribed tolerance~$\epsilon$,
providing a computable and KKT-consistent stopping rule.

\subsection{Convergence and Complexity Analysis}

Write the constrained problem as \(\min_x \Phi(x):=f(x)+\chi_F(x)\), where
\(\chi_F\) is the indicator of the closed convex set \(F\). Since
\(\Pi_F(z)=\argmin_x\{\chi_F(x)+\tfrac{1}{2}\|x-z\|_2^2\}\), NPGA is an
accelerated proximal-gradient method for \(\Phi=f+\chi_F\), with the proximal
map computed by the structured projector from Section~\ref{sec:npga}. The
theorem below specializes standard accelerated proximal-gradient theory to the
variables used throughout this paper, namely \(L_{\mathrm{eff}}\), \(\gamma\),
and the step-size sequence \(\{\alpha_k\}\). Part~(a) corresponds to
Algorithm~\ref{alg:npga} as written, whereas part~(b) records the classical
strongly convex constant-momentum specialization of the same framework; see
\citet{BeckTeboulle2009,Nesterov2004}.

\begin{theorem}[Convergence of NPGA]
\label{thm:npga-convergence}
Let \(f(x)=\|L_{\mathrm{eff}}^\top x\|_2^2+\gamma\|x\|_2^2\), where
\(L_{\mathrm{eff}}\) and \(\gamma\ge0\) specify the Baseline, Sketch, or STR
model. Define
\[
L_f = 2\bigl(\|L_{\mathrm{eff}}\|_2^2+\gamma\bigr),
\qquad
m_f = 2\bigl(\lambda_{\min}(L_{\mathrm{eff}}L_{\mathrm{eff}}^\top)+\gamma\bigr),
\]
and apply NPGA to \(\min_{x\in F}f(x)\) with exact projected steps
\(x^{k+1}=\Pi_F\bigl(y^k-\alpha_k\nabla f(y^k)\bigr)\). Let \(x^\star\) be any
minimizer of \(f\) over \(F\). Then:

\paragraph{(a) Convex case.}
Suppose \(m_f=0\). If NPGA uses the standard FISTA update
\(t_0=1\), \(t_{k+1}=(1+\sqrt{1+4t_k^2})/2\),
\(\beta_k=(t_k-1)/t_{k+1}\), together with either a fixed step
\(0<\alpha_k\equiv\alpha\le 1/L_f\) or the standard accelerated
proximal-gradient backtracking rule, then there exists
\(C_{\mathrm{cvx}}\ge0\), depending only on the initial point and the accepted
step sequence, such that
\[
f(x^k)-f(x^\star)
\;\le\;
\frac{C_{\mathrm{cvx}}}{(k+1)^2},
\qquad k\ge 1.
\]

\paragraph{(b) Strongly convex case.}
Suppose \(m_f>0\). If NPGA uses a fixed step \(0<\alpha\le 1/L_f\) and the
constant momentum
\[
\beta_k\equiv \beta:=\frac{1-\sqrt{\alpha m_f}}{1+\sqrt{\alpha m_f}},
\]
then there exists \(C_{\mathrm{lin}}\ge0\), depending only on the initial
point, such that
\[
f(x^k)-f(x^\star)
\;\le\;
C_{\mathrm{lin}}
\bigl(1-\sqrt{\alpha m_f}\bigr)^k,
\qquad k\ge 1.
\]

Thus the same exact-projection statement applies uniformly to the Baseline,
Sketch, and STR models through \(L_{\mathrm{eff}}\). When \(\gamma>0\), one has
\(m_f\ge 2\gamma>0\), so ridge stabilization is precisely what moves the
reduced model into the strongly convex linear-rate regime.
\end{theorem}
This theorem is not claimed as a new first-order proof in itself. Its role is
to express the standard smoothness and curvature constants directly through the
factor representation \(L_f=2(\|L_{\mathrm{eff}}\|_2^2+\gamma)\) and
\(m_f=2(\lambda_{\min}(L_{\mathrm{eff}}L_{\mathrm{eff}}^\top)+\gamma)\), so the
same solver statement covers the Baseline, Sketch, and STR models. In
particular, part~(a) is the FISTA-type guarantee for Algorithm~\ref{alg:npga}
itself, whereas part~(b) records the standard linear-rate specialization
available once \(m_f>0\). Existing
inexact variants such as \citet{SchmidtLeRouxBach2011} already treat
proximal-subproblem errors; what they do not directly encode is the
implementation-level stopping rule of our structured projector, which controls
the scalar dual search for \(x(\nu)=\Pi_\Delta(z+\nu\vmu)\).
Appendix~\ref{app:proofs-npga} summarizes this distinction in more detail.

\begin{remark}[Optimization error versus model error]
For a fixed reduced model \(M\), Theorem~\ref{thm:npga-convergence} gives the
solver error bound \(f_M(\hat x_M)-f_M(x_M^\star)\le \delta_k\), with
\(\delta_k=O(k^{-2})\) in part~(a) and \(\delta_k=O(\vartheta^k)\) in part~(b).
This is distinct from the modeling error caused by replacing the full factor
\(L\) with a sketched or STR surrogate. If
\(\sup_{x\in F}|f_M(x)-f_{\mathrm{full}}(x)|\le \Delta_M\), then
\(f_{\mathrm{full}}(\hat x_M)-f_{\mathrm{full}}(x_{\mathrm{full}}^\star)\le
\delta_k+2\Delta_M\). Section~\ref{sec:guarantees} provides exactly this
\(\Delta_M\)-type control for the Sketch and STR models.
\end{remark}

\subsection{GPU Parallelization Strategy}

The matrix--vector products in~\eqref{eq:npga-unified-grad}
are the dominant computational cost.
Because $L_{\mathrm{eff}}$ is stored in column-major GPU memory,
both $L_{\mathrm{eff}}^\top x$ and $L_{\mathrm{eff}}(L_{\mathrm{eff}}^\top x)$
can be computed using level-2 Basic Linear Algebra Subprograms (BLAS) kernels.
In practice, we fuse these two operations to minimize
memory traffic, where both multiplications are parallelized over
rows and columns using CUDA streams or batched general matrix--vector multiplication (GEMV) operations.
The projection $\Pi_F(v_k)$ is not coordinatewise separable, but it is built
from GPU-friendly bulk primitives: sorting, prefix sums, saxpy-style vector
updates, and dot products. In the active-constraint regime, each projection
call evaluates $x(\nu)=\Pi_\Delta(v_k+\nu\vmu)$ a small number of times inside
a scalar bracketing/bisection loop. This avoids large sparse factorizations and
keeps the projection step compatible with device-resident parallel kernels.

The key point here is structural rather than hardware-specific.
Once $L_{\mathrm{eff}}$ is resident in device memory, each NPGA iteration is dominated by two regular dense matrix--vector products and one projection step.
This makes the method well suited to BLAS- and sort-based GPU execution.
Actual runtime comparisons are reported later in Section~\ref{sec:experiments}.

\section{Theoretical Guarantees}
\label{sec:guarantees}

Section~\ref{sec:npga} controlled the optimization error of NPGA for any fixed
factor model \(L_{\mathrm{eff}}\). This section supplies the complementary
model-error bounds for the sketching and STR constructions, including
approximation, conditioning, and stability guarantees. Taken together, the two
sections yield the end-to-end accuracy decomposition used later in the
numerical experiments. Detailed proofs are given in
Appendices~\ref{app:proofs-embedding}--\ref{app:proofs-stability}.

\subsection{Subspace Embedding and Approximation Guarantees}
\label{sec:sketch}
We summarize the standard theorems that underpin our use of $\tilde L$ in place of $L$.

\begin{definition}[Subspace embedding]\label{def:embedding}
Let $U\in\mathbb{R}^{T\times r}$ have orthonormal columns spanning $\mathrm{Im}(L^\top)$, i.e., $U^\top U=I_r$ and $r=\rank(L)$. A random matrix $\Phi\in\mathbb{R}^{T\times s}$ is an $(\varepsilon,\delta)$-subspace embedding for $U$ if
\[
\Pr\!\left[\;\bigl\|(\Phi^\top U)^\top(\Phi^\top U) - I_r\bigr\|_2 \le \varepsilon\;\right]\ \ge\ 1-\delta.
\]
\end{definition}
\begin{remark}
The random matrix $\Phi$ in this definition does not depend on the specific choice of the orthonormal basis $U$ for $\mathrm{Im}(L^\top)$, but only on the subspace itself. This is because for any two orthonormal bases $U$ and $U'$ of $\mathrm{Im}(L^\top)$, there exists an orthogonal matrix $Q$ such that $U' = UQ$, and we have:
\[
(\Phi^\top U')^\top (\Phi^\top U') = Q^\top (\Phi^\top U)^\top (\Phi^\top U) Q.
\]
Since orthogonal transformations preserve the spectral norm ($\|Q^\top A Q\|_2 = \|A\|_2$ for any orthogonal $Q$), the embedding condition is invariant under the choice of basis. Therefore, $\Phi$ serves as a subspace embedding for the entire subspace $\mathrm{Im}(L^\top)$ regardless of the specific basis representation.    
\end{remark}
\begin{lemma}[Johnson--Lindenstrauss case {\cite{JohnsonLindenstrauss1984}}]
For $U\in\mathbb{R}^{T\times 1}$ with $\|U\|_2=1$, suppose
$s = O\bigl(\varepsilon^{-2}\log(1/\delta)\bigr)$. Then, with probability at
least $1-\delta$, $\bigl|\,\|\Phi^\top U\|_2^2 - 1\,\bigr| \le \varepsilon$.
\end{lemma}

\begin{theorem}[General subspace embedding {\cite{Sarlos2006,ClarksonWoodruff2017}}]\label{thm:embedding}
Let $U\in\mathbb{R}^{T\times r}$ with $U^\top U=I_r$, and suppose
$s = O\bigl((r+\log(1/\delta))/\varepsilon^2\bigr)$. Then, with probability at
least $1-\delta$, $\bigl\|(\Phi^\top U)^\top(\Phi^\top U)-I_r\bigr\|_2 \le
\varepsilon$.
\end{theorem}

\begin{corollary}[Geometry preservation on $\mathrm{Im}(L^\top)$]
\label{cor:geometry}
Let $\tilde L = L\Phi$. Assume that $\Phi$ is an $(\varepsilon,\delta)$-subspace embedding for $\mathrm{Im}(L^\top)$. Then, with probability at least $1-\delta$,
\[
(1-\varepsilon)\,\|L^\top x\|_2^2
\ \le\
\|\tilde L^\top x\|_2^2
\ \le\
(1+\varepsilon)\,\|L^\top x\|_2^2,
\qquad \forall\,x\in\mathbb{R}^n.
\]
\end{corollary}
A detailed proof is given in Appendix~\ref{app:proofs-embedding}.

\begin{remark}[Geometric interpretation]
Corollary~\ref{cor:geometry} states that the sketch preserves the norm of all vectors in the column space of $L^\top$ up to relative error $\varepsilon$. Equivalently, $\tilde\Sigma=\tilde L\tilde L^\top$ is a spectral approximation to $\Sigma=LL^\top$ when restricted to that subspace.
\end{remark}

This guarantees that replacing $L$ with $\tilde L$ preserves all quadratic 
forms associated with $\Sigma$, which is the key requirement for covariance 
approximation and downstream optimization.

\subsection{Error and Robustness Analysis}
We first control the spectral error between the original and sketched covariances. Such a bound is essential to justify using $\tilde\Sigma$ in place of 
$\Sigma$ for optimization, conditioning analysis, and regularization.

\begin{corollary}[Spectral error of the covariance approximation]
\label{cor:spectral}
If $\Phi$ is an $(\varepsilon,\delta)$-subspace embedding for $\mathrm{Im}(L^\top)$, then with probability at least $1-\delta$,
\[
\|\tilde\Sigma - \Sigma\|_2 \;\le\; \varepsilon\,\|\Sigma\|_2,
\qquad
\tilde\Sigma = \tilde L\tilde L^\top,\ \ \Sigma = LL^\top.
\]
In other words, $\tilde\Sigma$ \emph{spectrally approximates} $\Sigma$ within a
multiplicative factor of $(1\pm\varepsilon)$, preserving its eigenvalue
structure up to $O(\varepsilon)$ distortion.
\end{corollary}
A detailed proof is given in Appendix~\ref{app:proofs-embedding}.

Corollary~\ref{cor:spectral} further provides
quantitative insight into how the subspace embedding affects the eigenvalues of the covariance matrix.
\begin{corollary}[Positive definiteness preservation under subspace embedding]
\label{cor:pd-preserve}
Let $\Sigma=L L^\top$ be positive definite, with condition number
$\kappa(\Sigma)=\lambda_{\max}(\Sigma)/\lambda_{\min}(\Sigma)$, and let
$\tilde\Sigma=\tilde L \tilde L^\top$ be obtained by an
$(\varepsilon,\delta)$-subspace embedding for $\mathrm{Im}(L^\top)$. If
$\varepsilon < 1/\kappa(\Sigma)$, then with probability at least $1-\delta$,
$\tilde\Sigma$ remains positive definite.
\end{corollary}
A detailed proof is given in Appendix~\ref{app:proofs-embedding}.

\begin{theorem}[Optimal value error and solution sensitivity]
\label{thm:mv-robust}
Let
$v^\star=\min_{x\in F}x^\top\Sigma x$ and
$\tilde v^\star=\min_{x\in F}x^\top\tilde\Sigma x$.
If $\Phi$ is an $(\varepsilon,\delta)$-subspace embedding, then with probability at least $1-\delta$:
\[
\bigl|\tilde v^\star - v^\star\bigr|\ \le\ \varepsilon\,\|\Sigma\|_2.
\]
If, in addition, $\Sigma$ is positive definite and $x^\star,\tilde x^\star$ are minimizers for $\Sigma$ and $\tilde\Sigma$, respectively, then
\[
\|\tilde x^\star - x^\star\|_2 \;\le\; \varepsilon \kappa(\Sigma).
\]
\end{theorem}
A detailed proof is given in Appendix~\ref{app:proofs-stability}. The appendix
also records a restricted-strong-convexity extension of this stability result
for nearly singular covariance matrices.

\subsection{STR Conditioning and Approximation Guarantees}
\label{sec:STR}
Let \(\Sigma=Q\Lambda Q^\top\) with eigenvalues
\(\lambda_1\ge\cdots\ge\lambda_n\ge0\), and let
\(\Sigma_\ell:=Q_\ell\Lambda_\ell Q_\ell^\top\) denote the best rank-\(\ell\)
approximation of \(\Sigma\) in spectral norm. By the
Eckart--Young--Mirsky theorem \citep{HornJohnsonMatrixAnalysis},
\begin{equation}
\label{eq:rank-kbound}
\|\Sigma-\Sigma_\ell\|_2 = \lambda_{\ell+1}.
\end{equation}
The next theorem quantifies the spectral distortion of STR and gives a
sufficient condition for strict conditioning improvement.

\begin{theorem}[Conditioning improvement via STR]
\label{thm:L-STR}
Let $\tilde L=L\Phi$ be obtained by an $(\varepsilon,\delta)$-subspace embedding,
and let $\tilde L=U S V^\top$ be its thin SVD with
$S=\mathrm{diag}(\sigma_1,\dots,\sigma_r)$, $\sigma_1\ge\cdots\ge\sigma_r\ge0$.
Define the truncated factor
\[
\tilde L_\ell := U_\ell S_\ell V_\ell^\top,
\qquad
U_\ell := U_{[:,1:\ell]},\ \ 
V_\ell := V_{[:,1:\ell]},\ \ 
S_\ell=\mathrm{diag}(\sigma_1,\dots,\sigma_\ell).
\]
Let
\[
\widehat\Sigma := \tilde L_\ell \tilde L_\ell^\top + \gamma I_n,\qquad \gamma>0.
\]
Then, with probability at least $1-\delta$, the following hold:

\emph{(a) Spectral approximation bound.}
\begin{equation}
\label{eq:STR-spectral-bound}
\|\Sigma - \widehat\Sigma\|_2
\;\le\;
\|\Sigma - \Sigma_\ell\|_2 \;+\; 2\varepsilon\,\|\Sigma\|_2 \;+\; \gamma,
\end{equation}
where $\Sigma_\ell$ is the best rank-$\ell$ approximation of $\Sigma$ in spectral norm.

\emph{(b) Conditioning improvement on the covariance.}
Let 
\begin{equation}
\label{eq:gamma-kappa}
\gamma > 
\frac{(1+\varepsilon)\,\lambda_{\min}(\Sigma)\,\lambda_{\max}(\Sigma)}
{\lambda_{\max}(\Sigma)-(1+\varepsilon)\lambda_{\min}(\Sigma)}.
\end{equation}
then
\begin{equation}
\label{eq:STR-condSigma}
\kappa(\widehat\Sigma)\;<\;\kappa(\Sigma).
\end{equation}
\end{theorem}
A detailed proof is given in Appendix~\ref{app:proofs-stability}.
For any positive definite \(\Sigma\), adding \(\gamma I_n\) lowers the
condition number because
\((\lambda_{\max}+\gamma)/(\lambda_{\min}+\gamma)<\kappa(\Sigma)\) for every
\(\gamma>0\). The extra condition \eqref{eq:gamma-kappa} is needed only
because \(\widehat\Sigma\) is built from the perturbed low-rank proxy
\(\tilde L_\ell\tilde L_\ell^\top\). A practical discussion of choosing
\(\gamma\) is given in Appendix~\ref{app:str-discussion}.
\subsection{Stability of the Factorized Models}
\label{sec:lz}
We now analyze how approximating $L$ by sketching or STR
affects the optimal value and solution stability of the $L$--$z$ formulation.
Let $(x^\star,z^\star)$, $(\tilde x^\star,\tilde z^\star)$,
and $(\widehat x^\star,\widehat z^\star)$ be the minimizers of
the baseline, sketched, and STR-regularized problems, respectively:
\[
v^\star = \|L^\top x^\star\|_2^2,\quad
\tilde v^\star = \|\tilde L^\top \tilde x^\star\|_2^2,\quad
\widehat v^\star = \|\tilde L_\ell^\top \widehat x^\star\|_2^2 + \gamma\|\widehat x^\star\|_2^2.
\]

\paragraph{(1) Stability under subspace embedding.}
From Corollary~\ref{cor:spectral} and Theorem~\ref{thm:mv-robust},
if $\Phi$ is an $(\varepsilon,\delta)$-subspace embedding for $\mathrm{Im}(L^\top)$,
then with probability at least $1-\delta$,
\begin{equation}
\label{eq:Lz-stability}
\bigl|\tilde v^\star - v^\star\bigr|
\;\le\;
\varepsilon\,\|\Sigma\|_2,
\qquad
\|\tilde x^\star-x^\star\|_2
\;\le\;
\varepsilon\kappa(\Sigma).
\end{equation}
These bounds ensure that subspace embedding preserves both the
optimal value and the optimal solution up to $O(\varepsilon)$ accuracy.

\paragraph{(2) Stability under STR refinement.}
\begin{theorem}[Stability under STR approximation]
\label{thm:STR-stability}
Let $\Sigma = LL^\top$,
$
\tilde L = L\Phi \in \R^{n\times s},$ and $\widehat\Sigma = \tilde L_\ell \tilde L_\ell^\top + \gamma I_n,$
where $\Phi$ is an $(\varepsilon,\delta)$-subspace embedding,
$\tilde L_\ell$ is the top-$\ell$ truncation of $\tilde L$,
and $\gamma>0$ is the ridge parameter.
Then with probability at least $1-\delta$,
the following spectral-norm stability bound holds:
\begin{equation}
\label{eq:STR-stability-bound}
\boxed{\;
\|\widehat\Sigma - \Sigma\|_2
\;\le\;
2\varepsilon\,\|\Sigma\|_2
\;+\;
\frac{\tilde\lambda_{\ell+1}}{1-\varepsilon}
\;+\;
\gamma,
\;}
\end{equation}
where $\tilde\lambda_{\ell+1}$ is the $(\ell\!+\!1)$-th eigenvalue of
$\tilde\Sigma=\tilde L\tilde L^\top$.
Consequently, the optimal value and optimal solution of the STR model \eqref{eq:STR-model} satisfy
\begin{equation}
    \label{eq:STR-stability}
    \bigl|\widehat v^\star - v^\star\bigr|
\;\le\;
\|\widehat\Sigma - \Sigma\|_2,
\qquad
\|\widehat x^\star - x^\star\|_2
\;\le\;
\frac{\|\widehat\Sigma - \Sigma\|_2}{m},
\end{equation}
whenever $\Sigma\succeq mI_n$.
Equivalently,
\[
\|\widehat x^\star - x^\star\|_2
\;\le\;
\kappa(\Sigma)
\Big(
2\varepsilon
+
\frac{
\frac{\tilde\lambda_{\ell+1}}{1-\varepsilon} 
+\gamma}{\|\Sigma\|_2}
\Big).
\]
\end{theorem}
A detailed proof is given in Appendix~\ref{app:proofs-stability}.

\begin{remark}[Computational and numerical significance]
The $L$--$z$ factorization converts a dense quadratic form into a sparse
least-norm system that scales linearly in $n$. After sketching or STR, each
iteration uses only thin matrices ($n\times s$ or $n\times \ell$) and two
matrix--vector products. This reduces per-iteration complexity from $O(nT)$ to
$O(ns)$ (sketch) or $O(n\ell)$ (STR), while
\eqref{eq:Lz-stability}--\eqref{eq:STR-stability} guarantee that the resulting
solution and objective remain accurate.
\end{remark}

\section{Numerical Experiments}
\label{sec:experiments}

The numerical study asks three questions: how much approximation error do
sketching and STR introduce, how strongly does ridge improve conditioning, and
once the model is fixed, how competitive is NPGA against Gurobi on CPU and GPU?
We answer them through synthetic diagnostics, solver benchmarks, and
full-scale real-data tests.

\subsection{Experimental Design}
\label{sec:exp-design}

Section~\ref{sec:syn-approx-cond} addresses the first two questions on
synthetic data. It measures relative spectral error, full-objective gap, and
ridge-induced conditioning changes in order to connect the theory in
Section~\ref{sec:guarantees} to observable numerical behavior. These tests are
model-diagnostics experiments rather than throughput benchmarks.

Section~\ref{sec:syn-solver} then isolates the solver itself. On synthetic
instances we compare NPGA with Gurobi on the same training models, which lets
us separate first-order solver error from approximation error and examine when
GPU acceleration becomes materially useful.

Section~\ref{sec:real-bench} turns to the real 5-minute Chinese A-share data.
We first run subset benchmarks as a calibration step on moderate problem sizes,
where Gurobi is still comfortable. We then move to the full 5440-asset
instance, where the real computational value of the full, sketched, and STR
pipelines can be compared directly.

\subsection{Common Setup}
\label{sec:exp-setup}

Synthetic instances are generated from the controlled-spectrum model used
throughout the experiments. Unless stated otherwise, the factor singular values decay
geometrically, the synthetic return matrix has dimensions $n=600$ and
$T=2400$, and the effective rank is controlled by the retained-energy level
$\eta$ and the oversampling ratio $s/\ell$. Approximation errors are measured
against the empirical covariance $\Sigma=(1/T)RR^\top$. For a model $M$, solver
gaps are reported as
$$\mathrm{gap}_M(\hat x)=\frac{\max\{f_M(\hat x)-f_M(x_M^G),0\}}{\max\{|f_M(x_M^G)|,10^{-12}\}},$$
where $f_M$ is the regularized in-sample objective and $x_M^G$ is the
corresponding exact Gurobi solution. In
Section~\ref{sec:syn-approx-cond}, where the model itself changes, the reported
full-objective gap instead evaluates the approximate-model optimizer in the
unreduced objective $f_{\mathrm{full}}$; the ``Full-model gap'' column in
Table~\ref{tab:real-full} is defined analogously.

The real-data experiments use a cleaned RESSET 5-minute A--share panel from
January 2, 2020 to June 30, 2025. Missing vendor entries are filled with zero,
and assets with insufficient trading histories are removed, leaving a balanced
matrix with $n=5440$ stocks and $T=66412$ time periods. We split the sample at
January 1, 2024, which yields $48374$ training periods and $18038$ test
periods. Covariances are estimated from demeaned training returns, expected
returns from the raw training means, and the return threshold is set to the
$60$th percentile of the in-sample asset means. Test-period portfolio
statistics are annualized with $48\times 244=11712$ 5-minute intervals per
trading year.

All experiments are implemented in MATLAB R2023b, and exact QP benchmarks use
Gurobi 11.0.0. Runs were executed on a Dell Precision 3660 workstation running
Windows 11 with an Intel Core i9-12900K CPU, 128\,GB RAM, and an NVIDIA
GeForce RTX 4090 GPU with 24\,GB memory. The real-data panel comes from the
proprietary RESSET database. Runtime statistics are
reported as wall-clock times: when repeated timing runs are available, we
discard one warmup solve and report repeated-median timings, and in the
real-subset study each table entry averages three independently sampled
subsets. The approximation and conditioning diagnostics in
Section~\ref{sec:syn-approx-cond} are computed on the CPU, whereas GPU
acceleration is evaluated only in the solver benchmarks. We terminate NPGA when
$\|\Pi_F(x^k-\alpha_k\nabla f(x^k))-x^k\|_2$ falls below the prescribed
tolerance. For runtime benchmarks, timing runs use projection tolerance
$10^{-8}$ and recompute the KKT residual every five iterations, whereas
convergence traces retain per-iteration diagnostics. All reported NPGA
benchmarks use the fixed-step implementation based on the practical spectral
estimate described in Section~\ref{sec:npga}; the backtracking option is part
of the framework but is not separately benchmarked here.

\subsection{Synthetic Approximation and Conditioning}
\label{sec:syn-approx-cond}

\begin{figure}[t]
\centering
\begin{subfigure}[t]{0.8\linewidth}
\centering
\includegraphics[width=\linewidth]{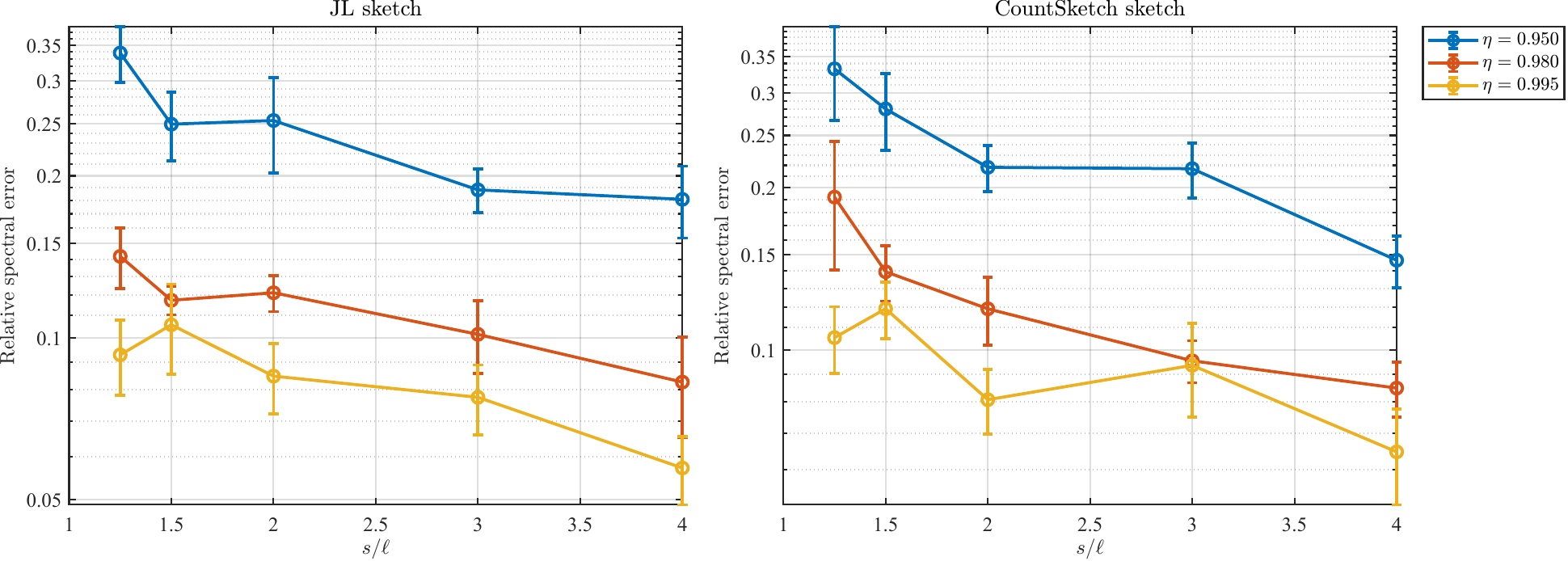}
\caption{Relative spectral error.}
\end{subfigure}

\medskip

\begin{subfigure}[t]{0.8\linewidth}
\centering
\includegraphics[width=\linewidth]{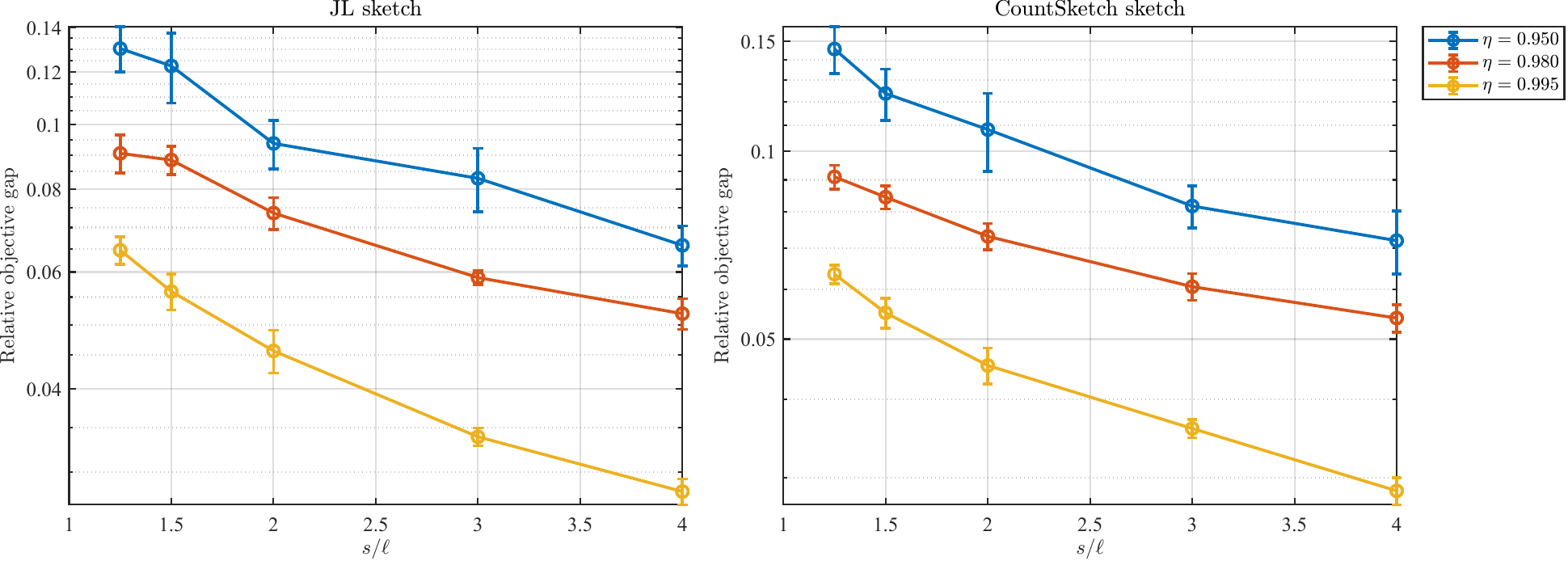}
\caption{Relative full-objective gap.}
\end{subfigure}
\caption{Approximation quality on synthetic instances as the sketch size grows.
Both metrics improve as $s/\ell$ increases, and JL and CountSketch behave
similarly once the retained energy level is fixed.}
\label{fig:syn-approx}
\end{figure}

Figure~\ref{fig:syn-approx} reports the relative spectral error and the
resulting training-objective gap as functions of $s/\ell$ for both JL and
CountSketch embeddings. These experiments are CPU diagnostics rather than GPU
benchmarks: they are designed to isolate approximation quality and ridge
stabilization before solver throughput enters the picture. The main message is
that the sketch dimension acts as a usable accuracy knob rather than as a
fragile heuristic. The trends are monotone in the expected direction:
larger $s/\ell$ and larger retained energy $\eta$ both reduce distortion. At
the representative setting $\eta=0.98$ and $s/\ell=2$, the two sketches are
nearly indistinguishable, with relative spectral error about $0.12$ and
relative objective gap about $7.3\%$; see Table~\ref{tab:syn-approx}. At the
more conservative setting $\eta=0.995$ and $s/\ell=4$, the mean objective gap
drops below $3\%$ for both sketch families. This behavior is consistent with
the perturbation bounds in Section~\ref{sec:guarantees}: once the ridge term
stabilizes the problem, covariance distortion translates into a controlled and
predictable loss in portfolio quality rather than an erratic degradation.

\begin{table}[ht]
\centering
\caption{Synthetic approximation quality at $\eta=0.98$ and $s/\ell=2$.}
\label{tab:syn-approx}
\begin{tabular}{lcc}
\toprule
Sketch & Rel.\ spec.\ err. & Rel.\ obj.\ gap\\
\midrule
JL & 0.1212 & 0.0736\\
CountSketch & 0.1192 & 0.0731\\
\bottomrule
\end{tabular}
\end{table}

Ridge regularization also improves conditioning markedly on a representative
CountSketch instance with $\eta=0.98$ and $s/\ell=2$. The sketched model starts
with condition number of order $10^8$ at $\gamma/\|\Sigma\|_2=10^{-8}$, but
this collapses to about $10^3$ once $\gamma/\|\Sigma\|_2$ reaches $10^{-3}$,
while $\lambda_{\min}$ increases by roughly five orders of magnitude. This is
exactly the numerical mechanism behind the strong-convexity regime analyzed in
Section~\ref{sec:npga}: the ridge term is not merely a technical device for the
proofs, but the practical lever that turns an unstable low-rank surrogate into
a well-conditioned optimization model. Detailed conditioning curves are
reported in Figure~\ref{fig:app-conditioning} in the appendix.

\subsection{Solver Accuracy and Scaling on Synthetic Problems}
\label{sec:syn-solver}

\begin{figure}[t]
\centering
\begin{subfigure}[t]{0.8\linewidth}
\centering
\includegraphics[width=\linewidth]{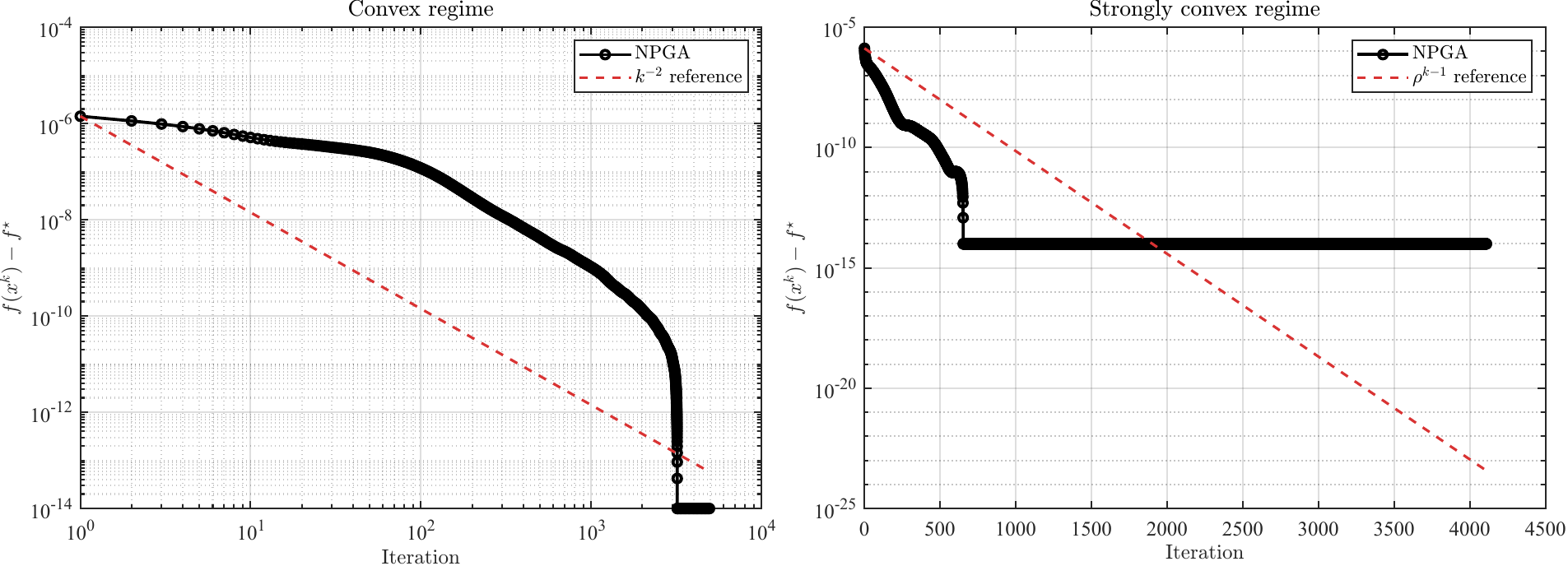}
\caption{Residual decay on log scales.}
\label{fig:syn-solver-resid}
\end{subfigure}

\medskip

\begin{subfigure}[t]{0.7\linewidth}
\centering
\includegraphics[width=\linewidth]{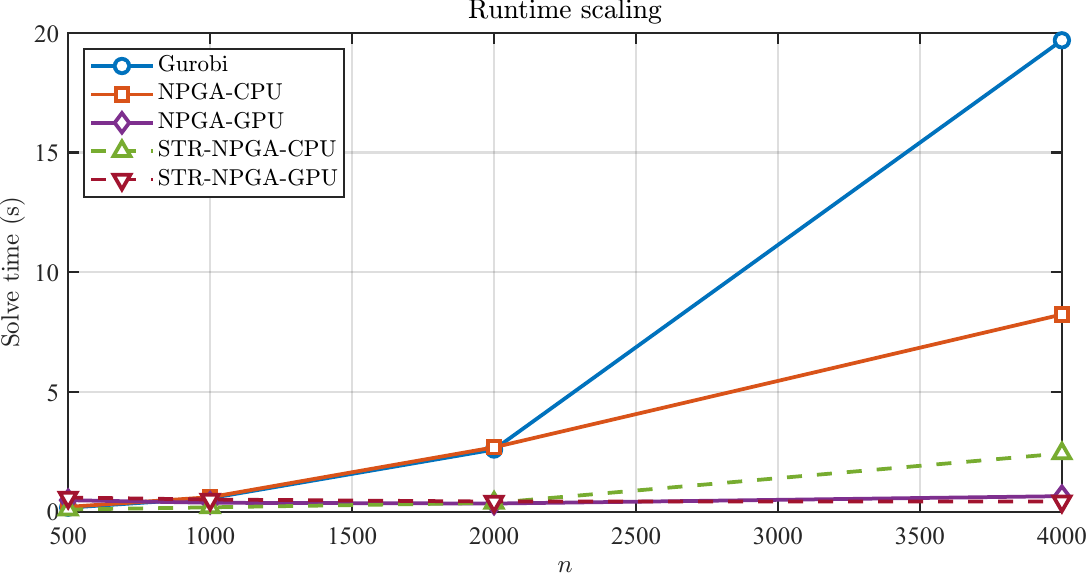}
\caption{Runtime scaling on log-log axes.}
\label{fig:syn-solver-runtime}
\end{subfigure}
\caption{Synthetic solver behavior. The convex case follows the
$k^{-2}$ benchmark, the ridge-stabilized case is close to the linear-rate
benchmark, and GPU acceleration becomes substantial on large unreduced
factors.}
\label{fig:syn-solver}
\end{figure}

Figure~\ref{fig:syn-solver-resid} plots the NPGA residual for the computed
iterates in a convex case and in a ridge-stabilized strongly convex case. The
observed curves track the reference $k^{-2}$ benchmark in the convex case and
the linear-rate benchmark in the ridge-stabilized strongly convex case. This
matches the two regimes highlighted by Theorem~\ref{thm:npga-convergence}:
part~(a) for the convex FISTA update, and part~(b) as the strong-convexity
benchmark made available once ridge stabilization yields \(m_f>0\). In the
reported runs we use the fixed-step FISTA implementation, so the linear curve
should be read as the benchmark associated with the strongly convex regime
rather than as a separate constant-momentum benchmark.
Figure~\ref{fig:syn-solver-runtime}
reports runtime scaling. GPU acceleration becomes decisive once the unreduced
factor is wide enough: at $n=4000$ and $T=16000$, NPGA-GPU solves the problem
in $0.64$ seconds versus $8.23$ seconds for NPGA-CPU, so even the unreduced
dense model becomes a subsecond-scale solve. On the compressed STR model,
STR-NPGA-GPU needs $0.42$ seconds versus $2.43$ seconds for STR-NPGA-CPU at
this largest size, although the smaller instances still show visible
device-side overhead. The practical implication is that GPU acceleration is
most valuable while the factor is still wide enough for dense linear algebra
to dominate the iteration cost.

Table~\ref{tab:syn-solver} isolates solver accuracy from approximation error by
comparing NPGA and Gurobi on the same training models. Once the problem size is
large enough, NPGA-GPU overtakes Gurobi decisively: at $n=2000$, Gurobi takes
$2.59$ seconds, while NPGA-GPU needs $0.33$ seconds with relative objective gap
$2.3\times 10^{-3}$; at $n=4000$, the gap widens to $19.69$ seconds versus
$0.64$ seconds. The same-model gap remains below $1.1\times 10^{-2}$
throughout the table, so the optimization error stays small relative to the
modeling error introduced later by sketching or truncation. STR-NPGA-GPU is
slightly slower than NPGA-GPU at $n=500,1000,2000$, but becomes faster at
$n=4000$ because the reduced factor finally compensates for the extra
projection overhead. These numbers show that the first-order solver is not the
weak link of the pipeline: once the model is fixed, NPGA reaches essentially
the same objective level as a commercial QP solver on the dense
full model at substantially lower cost.

\begin{table}[ht]
\centering
\caption{Synthetic solver benchmark against Gurobi on the same training problems.}
\label{tab:syn-solver}
\begin{tabular}{crrrrr}
\toprule
$n$ & Gurobi (s) & NPGA-GPU (s) & NPGA gap & STR-NPGA-GPU (s) & STR gap\\
\midrule
500 & 0.16 & 0.47 & 0.0006 & 0.59 & 0.0014\\
1000 & 0.56 & 0.37 & 0.0012 & 0.50 & 0.0025\\
2000 & 2.59 & 0.33 & 0.0023 & 0.42 & 0.0047\\
4000 & 19.69 & 0.64 & 0.0065 & 0.42 & 0.0108\\
\bottomrule
\end{tabular}
\end{table}

\subsection{Real-Data Benchmarks}
\label{sec:real-bench}

We first compare against Gurobi on random asset subsets drawn from the
5-minute Chinese A-share return data set. This subset experiment is only a
calibration step. Its role is to verify that NPGA tracks the exact Gurobi
solution closely on moderate instances where a general-purpose QP solver is
still comfortable on the dense full model. Each entry averages over
three random subsets.
In Table~\ref{tab:real-subset}, ``NPGA gap'' and ``STR gap'' both denote the
relative increase of the full-model in-sample objective above the Gurobi
benchmark on the same subset. The columns ``STR annualized return'' and
``STR annualized volatility'' report the out-of-sample annualized return and volatility of
the compressed STR portfolio; we only show these out-of-sample statistics for
STR because the main purpose of this subset test is to quantify the
approximation penalty of compression.
The table conveys two points. First, NPGA-CPU is numerically reliable on the
unreduced full model: its in-sample objective stays within $0.25\%$ to
$0.53\%$ of Gurobi, so the optimization error is already negligible on real
covariance structure even before GPU acceleration enters. Although on these
moderate subsets Gurobi remains faster, the attained objective values are
essentially the same. Second, STR-NPGA-CPU already cuts solve time below
Gurobi on every subset size, but it does so by accepting a visible
approximation error, with in-sample gaps between $6.7\%$ and $13.5\%$. Out of
sample, however, the degradation is milder: the STR portfolios have annualized
test volatility of $14.66\%$, $14.26\%$, and $12.79\%$ for
$n=300,600,1000$, which remains close to the corresponding Gurobi references.
The subset study therefore cleanly separates the two error sources: solver
error stays small, whereas the visible loss comes from model compression.

\begin{table}[ht]
\centering
\caption{Subset benchmarks on the real 5-minute A-share data.}
\label{tab:real-subset}
\resizebox{\linewidth}{!}{%
\begin{tabular}{crrrrrrr}
\toprule
$n$ & Gurobi (s) & NPGA-CPU (s) & NPGA gap & STR-NPGA-CPU (s) & STR gap & STR annualized return (\%) & STR annualized volatility (\%)\\
\midrule
300 & 0.22 & 1.28 & 0.0025 & 0.03 & 0.1348 & 36.15 & 14.66\\
600 & 0.62 & 2.10 & 0.0053 & 0.07 & 0.0937 & 39.30 & 14.26\\
1000 & 1.78 & 4.22 & 0.0045 & 0.12 & 0.0675 & 39.99 & 12.79\\
\bottomrule
\end{tabular}}
\end{table}

The full-scale benchmark is the decisive test because it is the setting that
motivates the paper. Here the training factor has $48374$ columns, whereas the
STR surrogate keeps only $\ell=474$ directions and uses a sketch of size
$s=948$ before truncation, so the final STR factor is more than $100\times$
smaller than the original one. Table~\ref{tab:real-full} now conveys three
complementary messages. First, on exactly the same $5440$-asset training
problem, NPGA-GPU solves the unreduced full model in $2.80$ seconds, versus
$64.84$ seconds for Gurobi and $26.77$ seconds for NPGA-CPU, while keeping the
full-model gap at $0.0080$ and essentially matching the Gurobi out-of-sample
portfolio statistics. This is the strongest evidence that the NPGA solver
itself is both accurate and fast on the original dense QP: the full model is
not merely solvable in principle, but operationally practical at this scale.

Second, the compressed GPU implementation improves materially once the runtime
benchmark uses the lighter projection configuration described above.
Sketch-NPGA-GPU drops to $3.68$ seconds end to end, and STR-NPGA-GPU drops to
$3.31$ seconds. These times indicate that, after compression, the dominant GPU
cost is not the matrix--vector products. It comes mainly from repeated
projections and synchronization around the projector.

Third, compression still changes the computational regime. STR-NPGA-CPU remains
the fastest end-to-end pipeline on this instance at $1.99$ seconds including
the sketch-and-truncate stage, and even the tuned compressed GPU rows remain
slower than full-model NPGA-GPU. Once the factor is reduced to $948$ or $474$
columns, the dense matvecs are no longer the dominant kernels; the remaining
cost is projection-dominated. Truncation nevertheless matters: reducing the
factor from $948$ to $474$ columns lowers the compressed GPU total time from
$3.68$ to $3.31$ seconds. The compressed surrogates remain approximate rather
than exact, with full-model objective gaps $0.2380$ for Sketch and $0.3203$ for
STR, but their out-of-sample return--volatility pairs remain close to the exact
reference. Additional convergence diagnostics in the appendix show the
same separation from two angles in Figure~\ref{fig:app-real-full-conv}: the
iteration counts stay of the same order,
whereas the solve times separate once projection cost dominates the compressed
models. The practical reading is therefore nuanced but favorable: if near-exact
full-model accuracy is the priority, NPGA-GPU is the strongest option; if the
goal is the fastest end-to-end reduced pipeline on this instance, STR-NPGA-CPU
still wins.

\begin{table}[ht]
\centering
\caption{Full-scale benchmark on the 5-minute A-share data with \(n=5440\)
assets and \(T=48374\) periods.} 
\label{tab:real-full}
\resizebox{\linewidth}{!}{%
\begin{tabular}{lrrrrrrr}
\toprule
Method & Factor cols & Build (s) & Solve (s) & Total (s) & Full-model gap & Annualized return (\%) & Annualized volatility (\%)\\
\midrule
Gurobi & 48374 & 0.00 & 64.84 & 64.84 & 0.0000 & 36.04 & 11.00\\
NPGA-CPU & 48374 & 0.00 & 26.77 & 26.77 & 0.0095 & 36.02 & 11.12\\
NPGA-GPU & 48374 & 0.00 & 2.80 & 2.80 & 0.0080 & 36.04 & 11.12\\
Sketch-NPGA-GPU & 948 & 1.15 & 2.53 & 3.68 & 0.2380 & 35.09 & 11.32\\
STR-NPGA-CPU & 474 & 1.15 & 0.83 & 1.99 & 0.3203 & 35.24 & 11.17\\
STR-NPGA-GPU & 474 & 1.15 & 2.16 & 3.31 & 0.3203 & 35.24 & 11.17\\
\bottomrule
\end{tabular}}
\end{table}

Taken together, the experiments support three conclusions. First, sketching
plus truncation provides a controllable accuracy--conditioning tradeoff rather
than an ad hoc reduction: moderate sketches keep the loss in solution quality
predictable, while ridge regularization makes the reduced models numerically
stable. Second, NPGA is accurate enough that solver error is not the dominant
source of loss, and NPGA-GPU already makes the unreduced dense full model
practical at the scale of the 5440-asset real-data benchmark. Third,
once the model is compressed aggressively, the computational bottleneck shifts
from dense linear algebra to projection and synchronization, which identifies
the next implementation target for further GPU gains.

\section{Conclusion and Future Work}
\label{sec:conclusion}

This paper presents a unified and scalable framework for large-scale
mean--variance portfolio optimization that integrates subspace embeddings,
$L$--$z$ factorization, theoretical guarantees, and a GPU-friendly
Nesterov-accelerated projected gradient solver. The main message is
threefold. First, sketching, truncation, and ridge regularization provide a
controllable approximation--conditioning tradeoff on ill-conditioned return
panels. Second, the theoretical analysis establishes explicit approximation,
conditioning, and stability guarantees for the sketching and STR models,
including $O(\varepsilon)$ bounds for covariance approximation, objective
value error, and solution perturbation under subspace embeddings. Third, NPGA
is accurate enough that the unreduced dense full model itself becomes
practical on modern GPUs: on the 5440-asset real-data benchmark, NPGA-GPU
solves the full model in $2.80$ seconds versus $64.84$ seconds for Gurobi
while preserving essentially the same out-of-sample portfolio statistics. The
compressed GPU experiments further show that, once the implementation is
tuned, the remaining bottleneck is projection and synchronization rather than
the factorized matrix--vector products. Future work includes stronger GPU
projection kernels, extensions to multi-factor or sparse portfolio models,
online or streaming settings, and dynamic state-estimation frameworks such as
the Yau--Yau filtering framework \citep{YauNiu2025,YauY00,YauY08}.

\section*{Acknowledgments}
The first author was supported by the National Natural Science Foundation of
China [Grant No.\ 42450242] and China's National Program of Overseas High-Level
Talent. Both authors also gratefully acknowledge institutional support from the
Beijing Institute of Mathematical Sciences and Applications (BIMSA).

\renewcommand{\maketitle}{}
\begin{APPENDICES}

\section{Relation of NPGA to Existing Accelerated Proximal-Gradient Theory}
\label{app:proofs-npga}

Write the constrained objective as
\[
\Phi(x):=f(x)+\chi_F(x),
\]
where \(\chi_F\) is the indicator of the closed convex feasible set \(F\). For
each extrapolated point \(y^k\), define the proximal subproblem objective
\[
\Psi_k(x):=
\langle \nabla f(y^k),x-y^k\rangle
+
\frac{1}{2\alpha_k}\|x-y^k\|_2^2
+
\chi_F(x).
\]
The exact projected NPGA step
\[
x^{k+1}=\Pi_F\bigl(y^k-\alpha_k\nabla f(y^k)\bigr)
\]
is exactly the proximal update for \(\Phi\) with stepsize \(\alpha_k\). In the
convex regime, Algorithm~\ref{alg:npga} therefore coincides with the standard
accelerated proximal-gradient / FISTA scheme applied to the composite
objective \(\Phi=f+\chi_F\), either with a fixed admissible step
\(0<\alpha\le 1/L_f\) or with the standard backtracking rule described in the
main text. The \(O(k^{-2})\) bound in Theorem~\ref{thm:npga-convergence}(a) is
exactly the classical FISTA rate specialized to \(h=\chi_F\); see Beck and
Teboulle (2009).

When \(m_f>0\), the model is strongly convex. If one switches from the
\(t_k\)-based momentum to the classical constant-momentum choice
\[
\beta=\frac{1-\sqrt{\alpha m_f}}{1+\sqrt{\alpha m_f}}
\]
with a fixed admissible step \(0<\alpha\le 1/L_f\), then one obtains the
linear-rate specialization recorded in
Theorem~\ref{thm:npga-convergence}(b). This is again a standard accelerated
proximal-gradient result; see Nesterov (2004) and Schmidt, Le Roux, and
Bach (2011).

The main paper does not claim a new accelerated-gradient proof here. Its
specific contribution is to express the smoothness and curvature conditions
directly through the factor representation
\[
L_f=2(\|L_{\mathrm{eff}}\|_2^2+\gamma),
\qquad
m_f=2(\lambda_{\min}(L_{\mathrm{eff}}L_{\mathrm{eff}}^\top)+\gamma),
\]
so that the same solver statement covers the Baseline, Sketch, and STR
quadratic models in one notation.

The point that is specific to the present paper is the structured projector.
General inexact accelerated proximal-gradient results, such as those of
Schmidt, Le Roux, and Bach (2011), are stated in terms of proximal-subproblem
errors. Concretely, they assume that the computed point \(x^{k+1}\) satisfies
\[
\Psi_k(x^{k+1})\le \inf_x \Psi_k(x)+\epsilon_k.
\]
Under such a condition, the standard inexact theory already gives perturbed
versions of the same FISTA-type rates, and likewise of the strongly convex
accelerated variants just mentioned.

Our implementation-level stopping rule is expressed differently. When the
return constraint is active, the projector is computed by scalar dual search in
\[
x(\nu)=\Pi_\Delta(z+\nu\vmu),
\]
and the inner loop is terminated once the one-dimensional residual
\[
\bigl|\vmu^\top x(\nu)-R_{\mathrm{target}}\bigr|
\]
falls below the prescribed scalar-search tolerance. This is a control on the
dual root-finding residual of the structured projector, not on the proximal
objective gap \(\Psi_k(x^{k+1})-\inf_x \Psi_k(x)\) that appears explicitly in
the generic inexact accelerated proximal-gradient theory.

Therefore, the only step not covered verbatim by the cited general results is
the projector-specific conversion from the scalar-search tolerance to a
proximal-subproblem error bound. One would need an additional lemma of the form
\[
\Psi_k\bigl(x(\hat\nu_k)\bigr)-\inf_x \Psi_k(x)\le \epsilon_k(\tau_k),
\]
where \(\tau_k\) denotes the scalar-search tolerance used to stop the dual
search. Once such a bound is available, the generic inexact accelerated
proximal-gradient results apply without further modification.

A useful first observation is that \(\Psi_k\) is \((1/\alpha_k)\)-strongly
convex, so if
\(p^k=\argmin_x \Psi_k(x)\), then every feasible point \(x\) satisfies
\[
\Psi_k(x)-\Psi_k(p^k)\ge \frac{1}{2\alpha_k}\|x-p^k\|_2^2.
\]
Hence any projector-specific estimate that bounds the proximal-objective gap
also immediately controls the distance to the exact proximal point. This
reduction is the only algorithm-specific ingredient missing from the standard
references, and it is not needed for the reported experiments because the
numerical runs in the main paper operate in the near-exact projection regime.
The full bibliographic details for these standard references appear in the main
paper.

\section{Proofs for Embedding and Covariance Approximation}
\label{app:proofs-embedding}

\subsection*{Proof of Corollary~\ref{cor:geometry}}
Let $U\in\mathbb{R}^{T\times r}$ be an orthonormal basis for the column space of
$L^\top$, i.e., $U^\top U=I_r$. Since $L^\top x\in\mathrm{Im}(L^\top)$, there
exists a unique $z\in\mathbb{R}^r$ such that
\[
L^\top x = U z, \qquad z = U^\top L^\top x.
\]
If $\Phi$ is an $(\varepsilon,\delta)$-subspace embedding for $U$, then
Theorem~\ref{thm:embedding} shows that with probability at least $1-\delta$,
\[
\bigl\|(\Phi^\top U)^\top(\Phi^\top U) - I_r\bigr\|_2 \;\le\; \varepsilon.
\]
Hence
\[
|z^{\top} \left((\Phi^\top U)^\top(\Phi^\top U) - I_r\right) z|
\leq
\bigl\|(\Phi^\top U)^\top(\Phi^\top U) - I_r\bigr\|_2 \|z\|_2^2
\leq
\varepsilon \|z\|_2^2, \quad \forall z\in \R^r.
\]
It follows that
\[
|\|\Phi^\top U z\|_2^2 - \|z\|_2^2| \leq \varepsilon \|z\|_2^2,
\quad \forall z\in \R^r,
\]
thus
\[
(1-\varepsilon)\,\|z\|_2^2
\ \le\
\|\Phi^\top U z\|_2^2
\ \le\
(1+\varepsilon)\,\|z\|_2^2,
\qquad \forall\,z\in\mathbb{R}^r.
\]
Substituting $z = U^\top L^\top x$ gives
\[
\Phi^\top U z = \Phi^\top L^\top x = (L\Phi)^\top x = \tilde L^\top x,
\qquad
\|z\|_2^2 = \|U^\top L^\top x\|_2^2 = \|L^\top x\|_2^2.
\]
Hence, with probability at least $1-\delta$,
\[
(1-\varepsilon)\,\|L^\top x\|_2^2
\ \le\
\|\tilde L^\top x\|_2^2
\ \le\
(1+\varepsilon)\,\|L^\top x\|_2^2,
\qquad \forall\,x\in\mathbb{R}^n.
\]
\Halmos

\subsection*{Proof of Corollary~\ref{cor:spectral}}
Let $\Delta := \tilde\Sigma - \Sigma = \tilde L \tilde L^\top - L L^\top$.
By the variational characterization of the spectral norm,
\[
\|\Delta\|_2
= \sup_{\|y\|_2=1}\,\bigl|\,y^\top \Delta\, y\,\bigr|
= \sup_{\|y\|_2=1}\,\bigl|\,\|\tilde L^\top y\|_2^2 - \|L^\top y\|_2^2\,\bigr|.
\]
Assume that $\Phi$ is an $(\varepsilon,\delta)$-subspace embedding for
$\mathrm{Im}(L^\top)$. By Corollary~\ref{cor:geometry}, with probability at
least $1-\delta$ we have for all $x\in\mathbb{R}^n$,
\[
(1-\varepsilon)\,\|L^\top x\|_2^2
\;\le\;
\|\tilde L^\top x\|_2^2
\;\le\;
(1+\varepsilon)\,\|L^\top x\|_2^2.
\]
Applying this inequality with $x=y$ for any unit vector $y$ yields
\[
\bigl|\,\|\tilde L^\top y\|_2^2 - \|L^\top y\|_2^2\,\bigr|
\;\le\;
\varepsilon\,\|L^\top y\|_2^2
\;=\;
\varepsilon\,y^\top(L L^\top)y
\;=\;
\varepsilon\,y^\top\Sigma\,y.
\]
Taking the supremum over all $\|y\|_2=1$ gives
\[
\|\Delta\|_2
\;=\;
\sup_{\|y\|_2=1}\bigl|\,y^\top \Delta\,y\,\bigr|
\;\le\;
\varepsilon\,\sup_{\|y\|_2=1}y^\top\Sigma\,y
\;=\;
\varepsilon\,\|\Sigma\|_2.
\]
Therefore, with probability at least $1-\delta$,
\[
\|\tilde\Sigma - \Sigma\|_2 \;\le\; \varepsilon\,\|\Sigma\|_2.
\]
Equivalently,
\[
(1-\varepsilon)\Sigma \preceq \tilde\Sigma \preceq (1+\varepsilon)\Sigma.
\]
\Halmos

\subsection*{Proof of Corollary~\ref{cor:pd-preserve}}
By Weyl's inequality and Corollary~\ref{cor:spectral}, for all $i=1,\ldots,n$,
with probability at least $1-\delta$,
\[
\bigl|\lambda_i(\tilde\Sigma)-\lambda_i(\Sigma)\bigr|
\;\le\;
\|\tilde\Sigma-\Sigma\|_2
\;\le\;
\varepsilon\,\|\Sigma\|_2
= \varepsilon\,\lambda_{\max}(\Sigma).
\]
Consequently,
\[
\lambda_{\min}(\tilde\Sigma)
\;\ge\;
\lambda_{\min}(\Sigma) -\varepsilon\,\lambda_{\max}(\Sigma),
\qquad
\lambda_{\max}(\tilde\Sigma)
\;\le\;
\lambda_{\max}(\Sigma) + \varepsilon\,\lambda_{\max}(\Sigma).
\]
If $\Sigma$ is positive definite and
$\varepsilon<\lambda_{\min}(\Sigma)/\lambda_{\max}(\Sigma)=1/\kappa(\Sigma)$,
then $\lambda_{\min}(\tilde\Sigma)>0$, ensuring that the sketched covariance
$\tilde\Sigma$ remains positive definite.
\Halmos

\section{Proofs for Robustness and Stability}
\label{app:proofs-stability}

\subsection*{Proof of Theorem~\ref{thm:mv-robust}}
Let $\Delta := \tilde\Sigma-\Sigma$. By Corollary~\ref{cor:spectral}, if $\Phi$
is an $(\varepsilon,\delta)$-subspace embedding for $\mathrm{Im}(L^\top)$, then
with probability at least $1-\delta$,
\begin{equation}
\label{eq:spectral-bound}
\|\Delta\|_2 \;=\; \|\tilde\Sigma-\Sigma\|_2 \;\le\; \varepsilon\,\|\Sigma\|_2 .
\end{equation}

\paragraph{(i) Optimal value error.}
Fix any $x\in F$. Since $x\ge 0$ and $\vone^\top x=1$, we have
$\|x\|_2\le \|x\|_1=1$. Then
\[
\bigl|\,x^\top \tilde\Sigma x - x^\top \Sigma x\,\bigr|
\;=\;\bigl|\,x^\top \Delta x\,\bigr|
\;\le\; \|\Delta\|_2\,\|x\|_2^2
\;\le\; \varepsilon\,\|\Sigma\|_2,
\]
where we used \eqref{eq:spectral-bound}. Let
$x^\star\in\argmin_{x\in F} x^\top \Sigma x$ and
$\tilde x^\star\in\argmin_{x\in F} x^\top \tilde\Sigma x$. Evaluating the
above bound at $x^\star$ and $\tilde x^\star$, we obtain
\[
\tilde v^\star
\,=\, \tilde x^{\star\top}\tilde\Sigma \tilde x^\star
\,\le\, \tilde x^{\star\top}\Sigma \tilde x^\star + \varepsilon\|\Sigma\|_2
\,\le\, x^{\star\top}\Sigma x^\star + \varepsilon\|\Sigma\|_2
\,=\, v^\star + \varepsilon\|\Sigma\|_2,
\]
and similarly,
\[
v^\star
\,=\, x^{\star\top}\Sigma x^\star
\,\le\, x^{\star\top}\tilde\Sigma x^\star + \varepsilon\|\Sigma\|_2
\,\le\, \tilde x^{\star\top}\tilde\Sigma \tilde x^\star + \varepsilon\|\Sigma\|_2
\,=\, \tilde v^\star + \varepsilon\|\Sigma\|_2.
\]
Combining the two displays yields
\[
\bigl|\tilde v^\star - v^\star\bigr| \;\le\; \varepsilon\,\|\Sigma\|_2.
\]

\paragraph{(ii) Solution sensitivity under strong convexity.}
Assume $\Sigma \succeq m I$ with $m:=\lambda_{\min}(\Sigma)>0$. Define
\[
f(x):=x^\top \Sigma x,
\qquad
\tilde f(x):=x^\top \tilde\Sigma x \;=\; f(x)+x^\top \Delta x .
\]
Then $\nabla f(x)=2\Sigma x$ and
$\nabla \tilde f(x)=2\tilde\Sigma x = \nabla f(x)+2\Delta x$.
By first-order optimality for convex optimization over the closed convex set
$F$,
\begin{equation}
\label{eq:VI}
\nabla f(x^\star)^\top (x - x^\star) \;\ge\; 0,
\qquad
\nabla \tilde f(\tilde x^\star)^\top (x - \tilde x^\star) \;\ge\; 0,
\quad \forall x\in F.
\end{equation}
Choose $x=\tilde x^\star$ in the first inequality of \eqref{eq:VI} and
$x=x^\star$ in the second, and add the two:
\begin{align*}
0
&\le \nabla f(x^\star)^\top(\tilde x^\star - x^\star)
 + \nabla \tilde f(\tilde x^\star)^\top (x^\star - \tilde x^\star) \\
&= \bigl(\nabla f(x^\star)-\nabla f(\tilde x^\star)\bigr)^\top(\tilde x^\star - x^\star)
   - \bigl(\nabla \tilde f(\tilde x^\star)-\nabla f(\tilde x^\star)\bigr)^\top(\tilde x^\star - x^\star) \\
&= \bigl(\nabla f(x^\star)-\nabla f(\tilde x^\star)\bigr)^\top(\tilde x^\star - x^\star)
   - \bigl(2\Delta \tilde x^\star\bigr)^\top(\tilde x^\star - x^\star).
\end{align*}
Rearranging gives
\begin{equation}
\label{eq:key-ineq}
\bigl(\nabla f(\tilde x^\star)-\nabla f(x^\star)\bigr)^\top(\tilde x^\star - x^\star)
\;\le\;
-2\,(\Delta \tilde x^\star)^\top(\tilde x^\star - x^\star).
\end{equation}
Applying Cauchy--Schwarz and \eqref{eq:spectral-bound} to the right-hand side
of \eqref{eq:key-ineq},
\begin{equation}
\label{eq:key-ineq2}
2\,\bigl|(\Delta \tilde x^\star)^\top(\tilde x^\star - x^\star)\bigr|
\;\le\;
2\,\|\Delta\|_2\,\|\tilde x^\star\|_2\,\|\tilde x^\star - x^\star\|_2
\;\le\;
2\,\varepsilon\,\|\Sigma\|_2\,\|\tilde x^\star - x^\star\|_2,
\end{equation}
where we used $\|\tilde x^\star\|_2 \le \|\tilde x^\star\|_1 = 1$ since
$\tilde x^\star \in F$. By $2m$-strong convexity of $f$,
\begin{equation}
\label{eq:strong-monotone}
\bigl(\nabla f(\tilde x^\star)-\nabla f(x^\star)\bigr)^\top(\tilde x^\star - x^\star)
\;\ge\; 2m\,\|\tilde x^\star - x^\star\|_2^2 .
\end{equation}
Combining \eqref{eq:key-ineq}, \eqref{eq:key-ineq2}, and
\eqref{eq:strong-monotone} yields
\[
2m\,\|\tilde x^\star - x^\star\|_2^2
\;\le\;
2\,\varepsilon\,\|\Sigma\|_2\,\|\tilde x^\star - x^\star\|_2,
\]
and hence, if $\tilde x^\star \neq x^\star$,
\[
\|\tilde x^\star - x^\star\|_2
\;\le\;
\frac{\varepsilon\,\|\Sigma\|_2}{m}
= \frac{\varepsilon \lambda_{\max}(\Sigma)}{\lambda_{\min}(\Sigma)}
= \varepsilon\kappa(\Sigma).
\]
The bound is trivial when $\tilde x^\star = x^\star$.
\Halmos

\subsection*{Restricted Strong Convexity Discussion}
When $\lambda_{\min}(\Sigma)$ is small or zero, the objective
$f(x)=x^\top\Sigma x$ is not globally strongly convex and the minimizer
of~\eqref{eq:mv} may be nonunique along flat directions. Nevertheless, many
structured problems exhibit restricted strong convexity on a subset of
directions that are relevant to the feasible set or to the true solution.
Formally, $f$ is said to satisfy an $m$-restricted strong convexity condition
with respect to a cone $\mathcal{C}\subseteq\mathbb{R}^n$ if
\[
f(x+h) - f(x) - \nabla f(x)^\top h
\;\ge\;
\frac{m}{2}\,\|h\|_2^2,
\qquad \forall\,h\in\mathcal{C}.
\]
Under RSC with parameter $m>0$, the stability conclusion in
Theorem~\ref{thm:mv-robust} continues to hold within $\mathcal{C}$: for
perturbations $\tilde\Sigma=\Sigma+\Delta$ satisfying
$\|\Delta\|_2\le\varepsilon\|\Sigma\|_2$,
\[
\|\tilde x^\star - x^\star\|_2
\;\le\;
\frac{\varepsilon\,\|\Sigma\|_2}{m},
\qquad
\text{for } \tilde x^\star - x^\star \in \mathcal{C}.
\]

\subsection*{Proof of Theorem~\ref{thm:L-STR}}
(a) Write $\Sigma=LL^\top$, $\tilde L=L\Phi$, and $\tilde\Sigma=\tilde L\tilde L^\top$.
Let $(\tilde\Sigma)_\ell$ denote the best rank-$\ell$ approximation to $\tilde\Sigma$. By the triangle inequality,
\[
\|\Sigma-\widehat\Sigma\|_2
\ \le\
\underbrace{\|\Sigma-\tilde\Sigma\|_2}_{\text{(I)}}\;+\;
\underbrace{\|\tilde\Sigma-(\tilde\Sigma)_\ell\|_2}_{\text{(II)}}\;+\;
\underbrace{\|(\tilde\Sigma)_\ell-\widehat\Sigma\|_2}_{\text{(III)}}.
\]
\emph{Term (I).} Since $\Phi$ is an $(\varepsilon,\delta)$-subspace embedding for $\mathrm{Im}(L^\top)$,
Corollary~\ref{cor:spectral} implies, with probability at least $1-\delta$,
\begin{equation}
\label{eq:termI}
\|\Sigma-\tilde\Sigma\|_2 \ \le\ \varepsilon\,\|\Sigma\|_2.
\end{equation}
\emph{Term (II).} By optimality of $(\tilde\Sigma)_\ell$ as the best rank-$\ell$ approximation to $\tilde\Sigma$,
\[
\|\tilde\Sigma-(\tilde\Sigma)_\ell\|_2 \ \le\ \|\tilde\Sigma-\Sigma_\ell\|_2
\ \le\ \|\tilde\Sigma-\Sigma\|_2+\|\Sigma-\Sigma_\ell\|_2,
\]
whence, using \eqref{eq:termI},
\begin{equation}
\label{eq:termII}
\|\tilde\Sigma-(\tilde\Sigma)_\ell\|_2 \ \le\ \varepsilon\,\|\Sigma\|_2 + \|\Sigma-\Sigma_\ell\|_2.
\end{equation}
\emph{Term (III).} Since $(\tilde\Sigma)_\ell=\tilde L_\ell\tilde L_\ell^\top$ and
$\widehat\Sigma=\tilde L_\ell\tilde L_\ell^\top+\gamma I_n$, we have
\[
(\tilde\Sigma)_\ell-\widehat\Sigma \;=\; -\,\gamma I_n
\quad\Rightarrow\quad
\|(\tilde\Sigma)_\ell-\widehat\Sigma\|_2 \;=\; \gamma.
\]
Combining \eqref{eq:termI}--\eqref{eq:termII} with the above yields
\eqref{eq:STR-spectral-bound} with probability at least $1-\delta$.

(b) The eigenvalues of $\widehat\Sigma$ are
$\{\sigma_1^2+\gamma,\dots,\sigma_\ell^2+\gamma,\underbrace{\gamma,\dots,\gamma}_{n-\ell}\}$,
hence
\[
\kappa(\widehat\Sigma)=\frac{\lambda_{\max}(\widehat\Sigma)}{\lambda_{\min}(\widehat\Sigma)}
=\frac{\sigma_1^2+\gamma}{\gamma}.
\]
Ensuring $\kappa(\widehat\Sigma)<\kappa(\Sigma)$ is equivalent to
\[
\frac{\sigma_1^2+\gamma}{\gamma}
\;<\;
\frac{\lambda_{\max}(\Sigma)}{\lambda_{\min}(\Sigma)}
\quad\Longleftrightarrow\quad
\gamma \;>\;
\frac{\lambda_{\min}(\Sigma)\,\sigma_1^2}{\lambda_{\max}(\Sigma)-\lambda_{\min}(\Sigma)}.
\]
Under an $(\varepsilon,\delta)$-subspace embedding we have the spectral bound
$\sigma_1^2(\tilde L)\le (1+\varepsilon)\lambda_{\max}(\Sigma)$ (with probability $\ge 1-\delta$),
which yields the following simpler sufficient condition:
\begin{equation}
\label{eq:gamma-kappa-sufficient}
\boxed{\;
\gamma \;>\;
\frac{(1+\varepsilon)\,\lambda_{\max}(\Sigma)\,\lambda_{\min}(\Sigma)}
{\lambda_{\max}(\Sigma)-(1+\varepsilon)\lambda_{\min}(\Sigma)}
\;}
\end{equation}
Then
\[
\kappa(\widehat\Sigma)
=\frac{\sigma_1^2+\gamma}{\gamma}
\;< \kappa(\Sigma),
\]
where the strict inequality is ensured by \eqref{eq:gamma-kappa-sufficient} after rearrangement.
\Halmos

\subsection*{Proof of Theorem~\ref{thm:STR-stability}}
Combining \eqref{eq:STR-spectral-bound} and \eqref{eq:rank_lbound}, we obtain
the composite stability bound \eqref{eq:STR-stability-bound}
\[
\|\Sigma - \widehat\Sigma\|_2
\;\le\;
\|\Sigma - \Sigma_\ell\|_2 \;+\; 2\varepsilon\,\|\Sigma\|_2 \;+\; \gamma
\;\le\;
\frac{\tilde\lambda_{\ell+1}}{1-\varepsilon} \;+\; 2\varepsilon\,\|\Sigma\|_2 \;+\; \gamma.
\]
Applying Theorem~\ref{thm:mv-robust} with $\Delta=\widehat\Sigma-\Sigma$ gives
\[
|\widehat v^\star-v^\star| \le \|\Delta\|_2,
\qquad
\|\widehat x^\star-x^\star\|_2 \le \|\Delta\|_2/m,
\]
whenever $\Sigma\succeq mI_n$. Combining with
\eqref{eq:STR-stability-bound} gives the stated results.
\Halmos

\section{Additional Discussion on STR}
\label{app:str-discussion}

\subsection{Practical Selection of the Truncation Level \texorpdfstring{$\ell$}{ell}}
\label{subsec:choice-l}

The truncation parameter $\ell$ determines the spectral split between the
signal-dominant directions and the noise-dominated bulk of the covariance
matrix. Instead of an energy-based rule, which is ineffective for real-world
financial covariances due to the extremely long, slowly decaying bulk, our
STR construction adopts a geometry-driven truncation principle derived from
the sketched factor matrix $\tilde L$. Let
\[
\tilde\lambda_i := \sigma_i^2(\tilde L),\qquad
\lambda_i := \sigma_i^2(L),
\]
sorted nonincreasingly. Under an $(\varepsilon,\delta)$-subspace embedding,
the singular values satisfy the spectral sandwich
\[
(1-\varepsilon)\lambda_i \le \tilde\lambda_i \le (1+\varepsilon)\lambda_i,
\qquad i=1,\dots,r.
\]
Thus the behavior of $\{\tilde\lambda_i\}$ reflects the latent structure of
$\{\lambda_i\}$ up to $(1\pm\varepsilon)$ distortion.

Empirically, $\{\tilde\lambda_i\}$ exhibits a clear core--bulk--tail pattern:
the first few singular values are factor-dominant and decay quickly, the
middle of the spectrum forms a long plateau associated with high-dimensional
noise, and the tail drops sharply toward the numerical floor. We therefore
choose $\ell$ at the end of the dominant core region, using the practical
rule
\[
\boxed{
\ell = \max\left\{\,i:\ 
\frac{\tilde\lambda_i}{\tilde\lambda_1} \ge \tau,\quad
\frac{\tilde\lambda_{i+1}}{\tilde\lambda_i} \le \rho \right\},
}
\]
where $\tau$ controls the minimum signal-to-noise ratio and $\rho<1$
detects the spectral knee separating the core from the bulk. Typical values
$\tau\in[10^{-3},10^{-2}]$ and $\rho\in[0.8,0.95]$ reliably identify the
factor-dominant directions in our experiments. This criterion is stable under
sketching because relative gaps $\tilde\lambda_{i+1}/\tilde\lambda_i$ are
preserved up to $1\pm O(\varepsilon)$.

Once $\ell$ is selected, the Eckart--Young--Mirsky theorem and the
$(\varepsilon,\delta)$-subspace embedding yield the computable upper bound
with probability at least $1-\delta$,
\begin{equation}
\label{eq:rank_lbound}
\boxed{\|\Sigma - \Sigma_\ell\|_2 = \lambda_{\ell+1}
\le \frac{\tilde\lambda_{\ell+1}}{1-\varepsilon}.}
\end{equation}
Thus the truncation error of STR can be controlled directly through the
observable sketched singular values without relying on total-variance
criteria.

\subsection{Choice of the ridge parameter \texorpdfstring{$\gamma$}{gamma}}
Once \(\ell\) is fixed, the STR covariance
\(\widehat\Sigma=\tilde L_\ell\tilde L_\ell^\top+\gamma I_n\) has condition
number
\[
\kappa(\widehat\Sigma)=\frac{\sigma_1^2(\tilde L_\ell)+\gamma}{\gamma}.
\]
Hence, if a target condition number \(\kappa_{\rm tar}>1\) is prescribed, one
may choose
\[
\gamma=\frac{\sigma_1^2(\tilde L_\ell)}{\kappa_{\rm tar}-1},
\]
which guarantees \(\kappa(\widehat\Sigma)\le \kappa_{\rm tar}\). In the real
data experiments, values such as \(\kappa_{\rm tar}\in[10^2,10^3]\) already
produce a substantial conditioning improvement while keeping the covariance
perturbation moderate.

\begin{remark}[Coupling between $\ell$ and $\gamma$]
The truncation level $\ell$ determines the smallest retained singular value
$\sigma_\ell$, while the ridge parameter $\gamma$ controls how much the
spectral floor is lifted. Since
\[
\kappa(\widehat\Sigma)=\frac{\sigma_1^2+\gamma}{\gamma},
\]
a smaller $\sigma_\ell$ (corresponding to a larger $\ell$) requires a larger
$\gamma$ to achieve any prescribed conditioning target. Therefore, $\ell$ and
$\gamma$ are complementary parameters: truncation filters noise-dominated
directions, and ridge regularization stabilizes the remaining spectrum.
\end{remark}

\subsection{Relation between STR and sketching}
The STR procedure can be interpreted as a structured sketching scheme that
extends classical randomized embeddings by incorporating deterministic
spectral regularization. Combining the randomized sketch stage with the
subsequent truncation and ridge stages, STR can be viewed conceptually as a
single well-conditioned low-rank sketching framework of the original factor
$L$:
\[
\widehat L \;=\; \bigl[\, L\Phi V_\ell \;\;,\;\; \sqrt{\gamma}\,I_n \,\bigr]
\in \R^{n\times(\ell+n)},
\]
so that
\begin{align*}
\widehat L\,\widehat L^\top
&\;=\;
(L\Phi V_\ell)(L\Phi V_\ell)^\top + \gamma I_n \\
&\;=\;
U S \bigl(V^\top V_\ell V_\ell^\top V\bigr) S U^\top + \gamma I_n \\
&\;=\;
U_\ell S_\ell^2 U_\ell^\top + \gamma I_n
\;=\;
\tilde L_\ell\tilde L_\ell^\top + \gamma I_n
\;=\;
\widehat\Sigma.
\end{align*}

\begin{remark}[Effect of STR on the sketch dimension $s$]
In the basic subspace embedding theory, the sketch size
\[
s = O\!\Big(\frac{r+\ln(1/\delta)}{\varepsilon^2}\Big)
\]
ensures $(\varepsilon,\delta)$-spectral approximation for the
$r$-dimensional column space of $L^\top$. When STR is applied, the subsequent
truncation and ridge steps operate only on the leading $\ell$ singular
directions of $\tilde L$. Therefore, the required sketch dimension depends on
the effective rank $\ell$ rather than the full rank $r$:
\[
s = O\!\Big(\frac{\ell+\ln(1/\delta)}{\varepsilon^2}\Big),
\qquad \ell\le r.
\]
\end{remark}

\subsection{Computational complexity of STR}
In large-scale settings, the effective rank $\ell$ and the true rank $r$ are
typically much smaller than the embedding dimension $s$, which in turn is far
smaller than the ambient dimensions $n$ and $T$:
\[
\ell \;\le\; r \;\le\; s \;\ll\; \min\{n,T\}.
\]
Here $n$ denotes the asset dimension, $T$ the sample length,
$s$ the embedding size, $r=\mathrm{rank}(\tilde L)$ the numerical rank of the
sketched factor, and $\ell$ the retained truncation level.

The first step constructs $\tilde L = L\Phi$ with a sketching matrix
$\Phi\in\R^{T\times s}$, where $s\ll\min\{n,T\}$. Its cost depends on the
structure of $\Phi$:
\[
\begin{array}{lcl}
\text{Gaussian JL:}  & O(nTs) & \text{(dense projection)};\\[2pt]
\text{CountSketch:}  & O(\mathrm{nnz}(L))\;\text{(typically }O(nT)\text{)}
& \text{(sparse, hash-based)}.
\end{array}
\]
The memory cost is $O(ns)$ to store $\tilde L$ (plus $O(Ts)$ for a dense
$\Phi$).

Computing the thin SVD $\tilde L=USV^\top$ requires $O(ns^2)$ flops and
$O(ns)$ memory. Truncating to the top-$\ell$ components
$(U_\ell,S_\ell,V_\ell)$ introduces negligible extra cost. For large-scale
implementations, one may instead compute a partial SVD via randomized or
Lanczos iterations, which scales as $O(ns\ell)$.

Forming $\widehat\Sigma = \tilde L_\ell\tilde L_\ell^\top + \gamma I_n$ is
conceptually $O(n\ell^2)$, but in practice the covariance is never explicitly
constructed. Overall, the one-shot STR cost is
\[
O(\underbrace{nT\!\cdot\!\text{cost}(\Phi)}_{\text{sketching}}
\;+\;
\underbrace{ns^2}_{\text{SVD/truncation}}),
\]
where $\text{cost}(\Phi)\in\{s,1\}$ for Gaussian JL and CountSketch,
respectively.

\section{Appendix Figures}
\label{app:figures}

This section collects the schematic and diagnostic plots referenced from the
main paper.
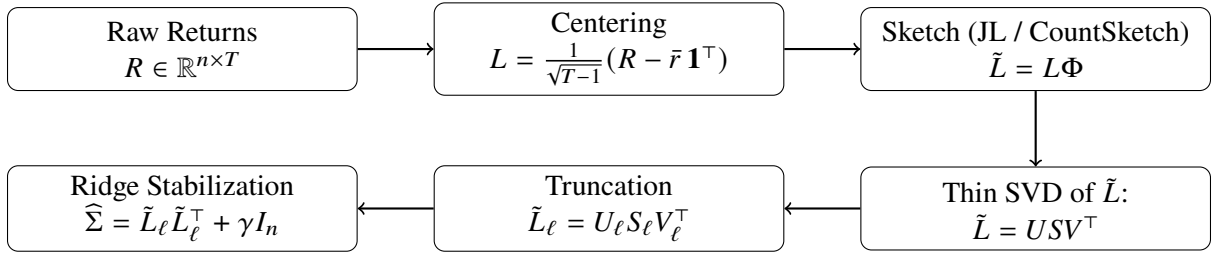
\begin{figure}[hb!]
\centering
\begin{tikzpicture}[
    node distance=1cm,
    box/.style={rectangle, rounded corners, draw, align=center,
                minimum width=4.6cm, minimum height=1.1cm},
    arrow/.style={->, thick}
]
\node[box] (R)         {Raw Returns \\ $R\in\mathbb{R}^{n\times T}$};
\node[box, right=of R] (Center)     {Centering \\ $L=\frac{1}{\sqrt{T-1}}(R-\bar r\,\mathbf{1}^\top)$};
\node[box, right=of Center] (Sketch) {Sketch (JL / CountSketch) \\ $\tilde L=L\Phi$};
\node[box, below=1cm of Sketch] (SVD)
    {Thin SVD of $\tilde L$: \\ $\tilde L = U S V^\top$};
\node[box, left=of SVD] (Trunc)
    {Truncation \\ $\tilde L_\ell = U_\ell S_\ell V_\ell^\top$};
\node[box, left=of Trunc] (Ridge)
    {Ridge Stabilization \\ $\widehat\Sigma=\tilde L_\ell\tilde L_\ell^\top+\gamma I_n$};
\draw[arrow] (R) -- (Center);
\draw[arrow] (Center) -- (Sketch);
\draw[arrow] (Sketch) -- (SVD);
\draw[arrow] (SVD) -- (Trunc);
\draw[arrow] (Trunc) -- (Ridge);
\end{tikzpicture}
\caption{The Sketch--Truncate--Ridge (STR) pipeline for constructing a stable covariance approximation.}
\label{fig:app-str-pipeline}
\end{figure}

\begin{figure}[hb!]
\centering
\includegraphics[width=0.88\linewidth]{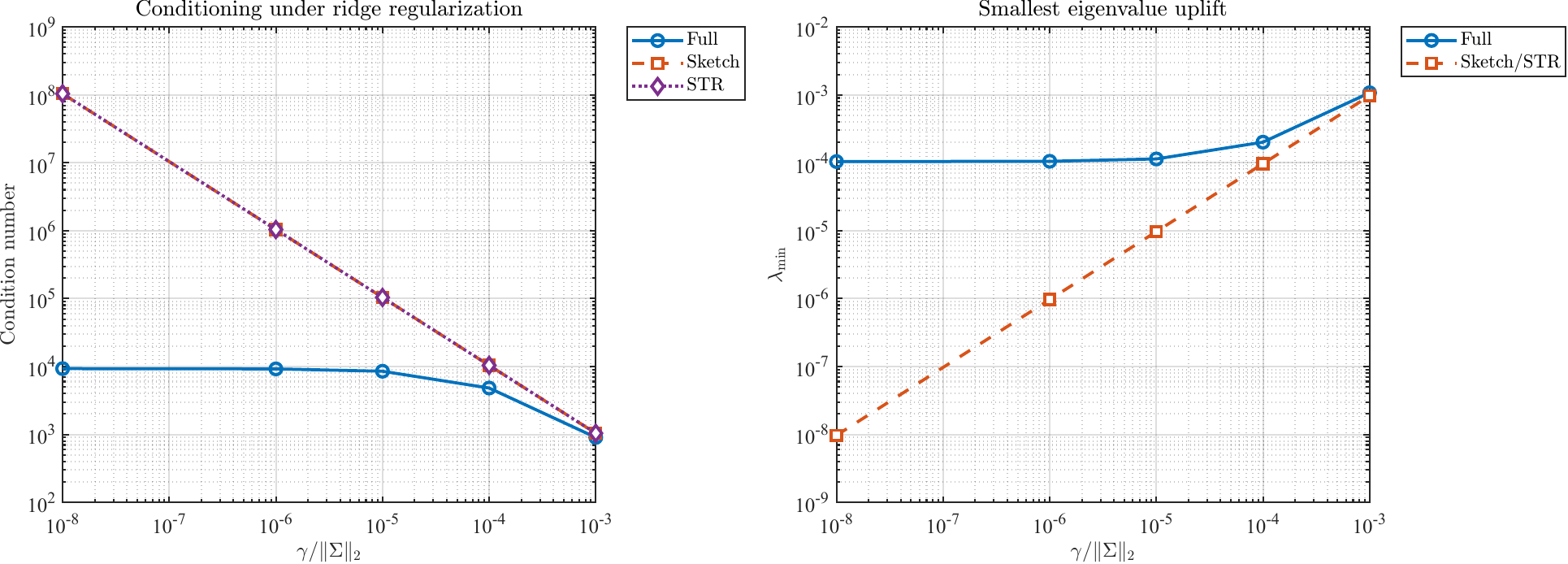}
\caption{Condition number and smallest-eigenvalue uplift under ridge regularization on a representative synthetic instance. The Sketch and STR curves coincide on this instance and therefore appear as a single trajectory.}
\label{fig:app-conditioning}
\end{figure}

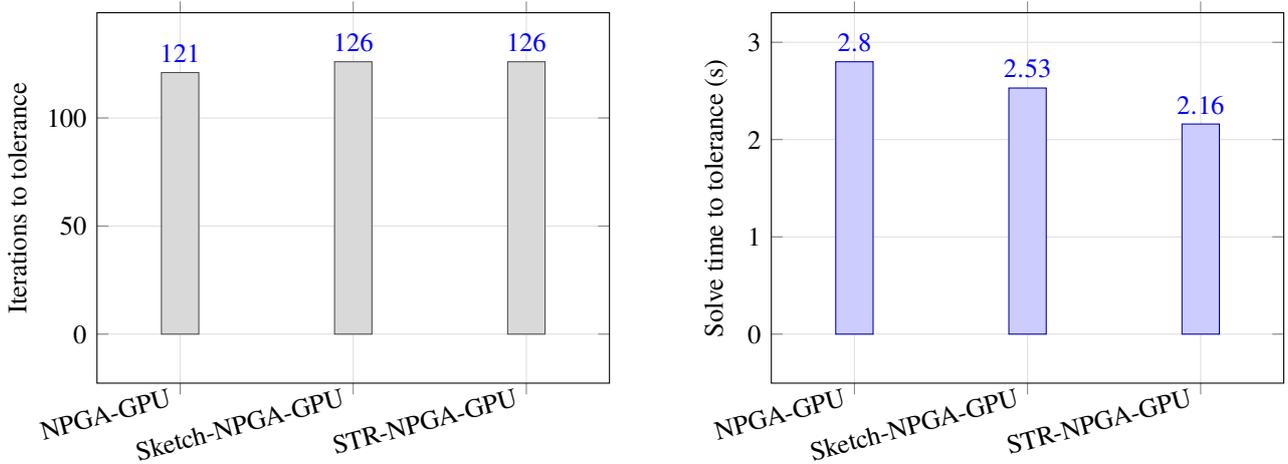
\begin{figure}[hb!]
\centering
\begin{subfigure}[t]{0.48\linewidth}
\centering
\begin{tikzpicture}
\begin{axis}[
width=\linewidth,
height=0.78\linewidth,
ybar,
ymin=0,
ylabel={Iterations to tolerance},
symbolic x coords={NPGA-GPU,Sketch-NPGA-GPU,STR-NPGA-GPU},
xtick=data,
xticklabels={NPGA-GPU,Sketch-NPGA-GPU,STR-NPGA-GPU},
grid=major,
major grid style={gray!25},
enlarge x limits=0.24,
tick label style={font=\small},
xticklabel style={font=\small, rotate=15, anchor=east},
label style={font=\small},
nodes near coords={\pgfmathprintnumber[fixed,precision=0]{\pgfplotspointmeta}},
every node near coord/.append style={font=\small, anchor=south},
point meta=y,
bar width=14pt,
enlarge y limits=0.18,
]
\addplot+[draw=black!70, fill=black!15]
coordinates {
(NPGA-GPU,121)
(Sketch-NPGA-GPU,126)
(STR-NPGA-GPU,126)
};
\end{axis}
\end{tikzpicture}
\caption{Outer iterations needed to reach the stopping tolerance.}
\end{subfigure}\hfill
\begin{subfigure}[t]{0.48\linewidth}
\centering
\begin{tikzpicture}
\begin{axis}[
width=\linewidth,
height=0.78\linewidth,
ybar,
ymin=0,
ylabel={Solve time to tolerance (s)},
symbolic x coords={NPGA-GPU,Sketch-NPGA-GPU,STR-NPGA-GPU},
xtick=data,
xticklabels={NPGA-GPU,Sketch-NPGA-GPU,STR-NPGA-GPU},
grid=major,
major grid style={gray!25},
enlarge x limits=0.24,
tick label style={font=\small},
xticklabel style={font=\small, rotate=15, anchor=east},
label style={font=\small},
nodes near coords={\pgfmathprintnumber[fixed,precision=2]{\pgfplotspointmeta}},
every node near coord/.append style={font=\small, anchor=south},
point meta=y,
bar width=14pt,
enlarge y limits=0.18,
]
\addplot+[draw=blue!70!black, fill=blue!20]
coordinates {
(NPGA-GPU,2.80)
(Sketch-NPGA-GPU,2.53)
(STR-NPGA-GPU,2.16)
};
\end{axis}
\end{tikzpicture}
\caption{Benchmark solve time to reach the stopping tolerance.}
\end{subfigure}
\caption{Convergence behavior of the GPU variants on the full 5440-asset real-data problem. The three GPU variants require nearly the same number of outer iterations, whereas their wall-clock times differ substantially.}
\label{fig:app-real-full-conv}
\end{figure}

\end{APPENDICES}

\end{document}